\documentclass[preprint, 11pt]{elsarticle}
\usepackage{graphicx}
\usepackage[utf8]{inputenc}
\usepackage{mdwlist}
\usepackage{amsfonts}
\usepackage{ulem}
\usepackage{amsmath}
\usepackage{amsthm}
\usepackage{amssymb}
\usepackage{mathbbol}
\usepackage{fancyhdr}
\usepackage{titlesec}
\usepackage{indentfirst}
\usepackage{booktabs}
\usepackage{verbatim}
\usepackage{color}
\usepackage{latexsym}
\usepackage{amscd}
\usepackage{esvect}
\usepackage{graphicx}
\usepackage{upgreek}
\usepackage{caption}   
\usepackage[page,header]{appendix}
\usepackage{titletoc}
\usepackage{mathrsfs}
\usepackage[numbers]{natbib}
\usepackage{float}
\usepackage{lipsum}
\usepackage[export]{adjustbox} 
\usepackage{subfiles}
\usepackage{caption}
\usepackage{subcaption}

\usepackage{algorithm}
\usepackage{algorithmicx}
\usepackage[noend]{algpseudocode}

\algnewcommand\algorithmicinput{\textbf{INPUT: }}
\algnewcommand\Input{\item[\algorithmicinput]}
\algnewcommand\algorithmicoutput{\textbf{OUTPUT: }}
\algnewcommand\Output{\item[\algorithmicoutput]}

\numberwithin{equation}{section}

\allowdisplaybreaks

\newcommand{\beq}{\begin{equation}}
\newcommand{\eeq}{\end{equation}}
\newcommand{\ben}{\begin{eqnarray}}
\newcommand{\een}{\end{eqnarray}}
\newcommand{\beno}{\begin{eqnarray*}}
\newcommand{\eeno}{\end{eqnarray*}}
\newcommand{\bz}{\mathbf{z}}
\def\veps{\varepsilon}

\def\bz {\boldsymbol{z}}
\def\e{\varepsilon}

\begin{document}

\begin{frontmatter}
\title{On a neural network approach for solving potential control problem of the semiclassical Schr\"odinger equation}
\author[1]{Yating Wang}
\author[2]{Liu Liu}
\address[1]{School of Mathematics and Statistics, Xi'an Jiaotong University.}
\address[2]{Department of Mathematics, Chinese University of Hong Kong.}

\begin{abstract}
    Robust control design for quantum systems is a challenging and key task for practical technology. In this work, we apply neural networks to learn the control problem for the semiclassical Schr\"odinger equation, where the control variable is the potential given by an external field that may contain uncertainties. Inspired by a relevant work \cite{NNLC}, we incorporate the sampling-based learning process into the training of networks, while combining with the fast time-splitting spectral method for the Schr\"odinger equation in the semi-classical regime. The numerical results have shown the efficiency and accuracy of our proposed deep learning approach. 
\end{abstract}

\end{frontmatter}

\section{Introduction}

Control of quantum phenomena has been an important scientific problem in the emerging quantum technology \cite{Survey}. The control of quantum electronic states in physical systems has a variety of applications such as quantum computers \cite{BBBV97}, control of photochemical processes \cite{BDMA06} and semiconductor lasers \cite{laser}. Detailed overviews of the quantum control field can be found in survey papers and monographs \cite{DP10, WM09}. 
One issue of the controllability theory \cite{Control-system} aims to assess the ability to steer a quantum system from an arbitrary initial state to a targeted final state, under the impact of a control field such as a potential function, given possibly noisy observation data. 

Uncertainty Quantification (UQ) has drawn many attentions over the past decade. In simulating physical systems, which are often modeled by differential equations, there are inevitably modeling errors, imprecise measurements of the initial data or background coefficients, which may bring uncertainties to the models. In this project, we study the semiclassical Schr\"odinger equation with external potential that may contain uncertainties, and is treated as the control variable. Let $\Omega$ be a bounded domain in $\mathbb R$, the Schr\"odinger equation in the semiclassical regime is described by a wave function 
$\psi: \mathcal{Q} \mapsto \mathbb{C}$, 
\begin{equation}
\label{Schro-eqn}
\left\{
\begin{array}{ll}
\displaystyle i \e \partial_t \psi^{\e} = - \frac{\e^2}{2} \Delta \psi^{\e} + V(x,\bz) \psi^{\e}, \qquad  (x, t) \in \mathcal{Q}\times (0,T),  \\[8pt]
\displaystyle \psi |_{t=0} = \psi_0(x), \qquad x\in\Omega \subset \mathbb{R}, 
\end{array}
\right.
\end{equation}
where $0<\e \ll 1$ is the scaled Planck constant describing the microscopic and macroscopic scale ratio. Here the solution $\psi = \psi(t,x, \bz)$ is the electron wave function with initial condition $\psi_0(x)$, the potential $V(x, \bz)\in L^{\infty}(\Omega\times I_{\bz})$ is the control variable that models the external field and is spatially dependent. Periodic boundary condition is assumed in our problem. 

The uncertainty is described by the random variable $\bz$, which lies in the random space $I_{\bz}$ with a probability measure $\pi(\bz)d \bz$.  
We introduce the notation for the expected value of $f(\bz)$ in the random variable $\mathbf z$, 
\begin{equation}
\langle f \rangle_{\pi(\mathbf z)} = \int f(\bz) \pi(\mathbf z) d \mathbf z. 
\end{equation}
The solution to the Schr\"odinger equation is a complex valued wave function, whose nonlinear transforms lead to probabilistic measures of the physical observables. The primary physical quantities of interests include position density,
\begin{equation}
\label{density1} 
n^{\e} = |\psi^{\e}|^2, 
\end{equation}
and current density
\begin{equation}
\label{density2}
J^{\e} = \e \operatorname{Im}\left(\overline{\psi^{\e}}\nabla\psi^{\e}\right) 
= \frac{1}{2i} \left( \overline{\psi^{\e}}\nabla \psi^{\e} - \psi^{\e}\nabla\overline{\psi^{\e}}\right). 
\end{equation}
At each fixed $\bz$, with $V$ being continuous and bounded, the Hamiltonian operator $H^{\e}$ defined by 
$$ H^{\e}\psi^{\e} = - \frac{\e^2}{2} \Delta \psi^{\e} + V(x,\bz) \psi^{\e} $$
maps functions in $H^2(\mathbb{R}^d)$ to $L^2(\mathbb{R}^d)$ and is self-adjoint. The operator 
$\frac{1}{i \e} H^{\e}$ generates a unitary, strongly continuous semi-group on $L^2(\mathbb{R}^d)$, which guarantees a unique solution of the Schr\"odinger equation \eqref{Schro-eqn} that lie in the space \cite{WBV10}: 
\begin{equation*}
W(0,T) := \left\{ \phi \in L^2((0,T); H_0^1(\Omega;\mathbb{C})) \Big| \frac{d\phi}{dt}
\in L^2((0,T); H^{-1}(\Omega;\mathbb{C})) \right\}. 
\end{equation*}

As a literature review, we mention that there has been several work \cite{BP02,LZ02} on boundary control for the Schr\"odinger equation \eqref{Schro-eqn}, where the observation is taken from the Dirichlet or Neumann boundary data. In some references such as \cite{Sampling-control}, the authors consider the quantum system with evolution of its state $|\psi(t)\rangle$ described by the Schr\"odinger equation $\dfrac{d}{dt}|\psi(t)\rangle = -i H(t)|\psi(t)\rangle$ with the initial condition $|\psi(0)\rangle = |\psi_0 \rangle$. 
The Hamiltonian $H(t)$ there corresponds to a time-dependent control variable that contains random parameters. Their goal is to drive the quantum ensemble from an initial state $|\psi_0\rangle$ to the target state $|\psi_{\text{target}}\rangle$, by employing a gradient-based learning method to optimize the control field. In \cite[Section 7.3]{BCS17}, the control problem of a charged particle in a well potential was formulated, where in their setting the potential field is time-dependent. We mention some other relevant work on stability estimates and semiclassical limit of inverse problem for the Schr\"odinger equation \cite{BF10, CL22, Eskin08, Lemm, WBV10}. 

{\textcolor{black}{We continue to mention several studies that are related to the inverse problems for the Schr\"odinger equation or other models.}}
For relevant inverse boundary value problems on this topic, there are existing iterative methods applied to the Helmholtz equation \cite{HM03}, where one starts with an initial guess of the boundary condition, then adjusts it iteratively by minimizing functionals such as error norms between the calculated data and measured data. This could be extremely time-consuming since at each iteration step, a forward problem needs to be solved. In the partial boundary data situation, there has been research on studying the linearized inverse problem of recovering potential function for the 
time-independent Schr\"odinger equation \cite{ZLX22}. Moreover, for inverse potential problems, well-posedness of the continuous regularized formulation was analyzed in both elliptic and parabolic problems, with conditional stability estimates and error analysis for the discrete scheme studied in \cite{CY08, Bangti}. 

{\color{black}{The desired control problem can be described as the following: To which extend can the wave solution 
$\psi^{\e}$ of \eqref{Schro-eqn} be perturbed by the control field--in our case the potential function $V$, in order to reach the desired target state at the final time $T$?}} The above question can be reformulated into an {\it optimal control} problem. At the final time $T$, given the target state $\psi_{{\it target}}$, let $V$ be approximated by a neural network parameterized by $\boldsymbol{\theta}$, and $\lambda>0$ be a regularization coefficient, we aim to solve the following minimization problem: 
{\color{black}
\begin{equation}
\label{Opt}
\left\{
\begin{array}{lll}
\displaystyle\min_{\boldsymbol{\theta}} J_{\lambda}(V(\boldsymbol{\theta})) = \min_{\boldsymbol{\theta}}   ||\psi^{\e}(x,T;\boldsymbol{\theta}) - \psi_{{\it target}}||_{L^2(\Omega)}^2 + \lambda\, || V(x;\boldsymbol{\theta}) ||_{L^2(\Omega)}^2, \\[6pt]
\text{such that   }\displaystyle \quad i \e \partial_t \psi^{\e}(x,t;\boldsymbol{\theta}) = - \frac{\e^2}{2} \Delta \psi^{\e}(x,t;\boldsymbol{\theta}) + V(x;\boldsymbol{\theta}) \psi^{\e}(x,t;\boldsymbol{\theta}), \\[6pt]
\qquad\qquad\quad \psi^{\e}(x,t=0;\boldsymbol{\theta}) = \psi_0(x). 
\end{array}
\right.
\end{equation}
if $V$ is a deterministic potential, and 
\begin{equation}
\label{Opt}
\left\{
\begin{array}{lll}
\displaystyle\min_{\boldsymbol{\theta}} J_{\lambda}(V(\boldsymbol{\theta})) = \min_{\boldsymbol{\theta}}   || \psi^{\e}(x,T;\boldsymbol{\theta},\bz) - \psi_{{\it target}}(\bz)||_{L^2(\Omega\times I_{\bz})}^2 + \lambda\, || V(x;\boldsymbol{\theta},\bz) ||_{L^2(\Omega\times I_{\bz})}^2 , \\[6pt]
\text{such that   }\displaystyle \quad i \e \partial_t \psi^{\e}(x,t;\boldsymbol{\theta},\bz) = - \frac{\e^2}{2} \Delta \psi^{\e}(x,t;\boldsymbol{\theta},\bz) + V(x;\boldsymbol{\theta},\bz) \psi^{\e}(x,t;\boldsymbol{\theta,\bz}), \\[6pt]
\qquad\qquad\quad \psi^{\e}(x,t=0;\boldsymbol{\theta},\bz) = \psi_0(x;\bz). 
\end{array}
\right.
\end{equation}
if the potential $V$ contains uncertainty and the random variable is $\bz$.
}

In each particular problem setting, discretized form of the above loss function will be presented. 
We now highlight the main contributions of our work: 
\begin{enumerate}
    \item We take advantage of the rising trend of machine learning and use neural networks to approximate the control variable considered as the potential field in the Schr\"odinger equation. Both deterministic and stochastic control functions are considered. A fully-connected neural network is used for the deterministic problem, and the DeepONet \cite{LJPZG21} is applied in the stochastic case.
    \item During the training process, the Schr\"odinger equation in the semiclassical regime is solved using the fast time-splitting spectral method to improve the computational efficiency and accuracy of our algorithm. 
    \item We study and compare both cases when the observation data is associated with or without noise, and propose different training strategies. For data without noise, the popular stochastic gradient descent (SGD) method is used. For noisy data, we consider a Bayesian framework and adopt the stochastic gradient Markov chain Monte Carlo (MCMC) approach to obtain robust learning results. 
\end{enumerate}

The rest of the paper is organized as follows. In Section \ref{sec:regularity}, we discuss the oscillatory behavior of solution to the semiclassical Schr\"odinger equation in the random variable and mention the numerical challenges even for the forward UQ problems. 
Our main methodology of using the learning-based technique to solve the optimization problem \eqref{Opt} will be proposed in Section \ref{sec:method}, with numerical scheme for the forward problem introduced in subsection \ref{subsec:forward} and several neural network approaches described in subsection \ref{subsec:net}. We conduct extensive numerical experiments for both the deterministic and stochastic potential control problems and present the results in Section \ref{sec:numerical}. Conclusion and future work will be addressed lastly. 

\section{Regularity of solution in the random space}\label{sec:regularity}

The semi-classical Schr\"odinger equation is a family of dispersive wave equations parameterized by $\e \ll 1$, it is well known that the wave equation propagates $O(\e)$ scaled oscillations in space and time. However, for UQ problems it is not obvious whether the small parameter $\e$ induces oscillations in the random variable $\mathbf z$. We conduct a regularity analysis of $\psi$ in the random space, which enables us to study the oscillatory behavior of solution in the random space. 

To investigate the regularity of the wave function in the $\bz$ variable, we check the following averaged norm 
\begin{equation}
\label{energy} ||\psi||_{\Gamma }:= \left (\int_{I_z}\int_{\mathbb R^3}\left| \psi(t, \mathbf x, \bz)\right|^2\,  d{\mathbf x}\pi(\bz)d{\bz} \right)^{1/2}. 
\end{equation}
First, observe that
$\forall\, \mathbf z \in I_{\mathbf z}$, 
\[
	\frac{\partial}{\partial t} \| \psi^{\e}\|^2_{L^2_{\mathbf x}}(t,\bz)=0, 
\]
thus 
\[
\frac{d}{dt} \| \psi^{\e}\|_{\Gamma}^2=0,
\]
which indicates the $\Gamma$-norm of the wave function $\psi^\e$ is conserved in time,
$ \psi^{\e}\|_{\Gamma} (t) =  \| \psi^{\e}_{\text{in}}\|_{\Gamma}\,. $

Below we show that $\psi^\e$ has $\e$-scaled oscillations in $\mathbf z$. As an example, we analyze first-order partial derivative of $\psi^\e$ in $z_1$ and denote $\psi^1= \psi^\e_{z_1}$ and $V^1=V_{z_1}$. By differentiating the semi-classical Schr\"odinger equation \eqref{Schro-eqn} with respect to $z_1$, one gets
\[
    i \e \psi^1_t =- \frac{\e^2}{2}\Delta_{\mathbf x} \psi^1 + V^1 \psi^\e + V \psi^1.
\]
Direct calculation leads to
\begin{align*}
\frac{d}{dt} \| \psi^{1}\|^2_{\Gamma} &= \int \bigl (\psi^1_t \bar \psi^1 + \psi^1 \bar \psi^1_t \bigr) \pi d \mathbf x d \mathbf z \\
& = \int \bigl(\frac{1}{i\e} V^1 \psi^\e \bar \psi^1  - \frac{1}{i\e} V^1 \psi^1 \bar \psi^\e \bigr) \pi d \mathbf x d \mathbf z \\
& \le \frac{2}{\e} \| \psi^{1}\|_{\Gamma}\, \| V^1 \psi^{\e}\|_{\Gamma}\,,
\end{align*}
where we use the Cauchy-Schwarz inequality and Jensen inequality in the last step, namely
\begin{align*}
&\int V^{1}\psi^{\e}\bar\psi^{1} dx \leq \left( \int (V^{1}\psi^{\e})^2 dx \right)^{1/2} \left( \int (\bar\psi^{1})^2 dx\right)^{1/2}, \\
&\int\int V^{1}\psi^{\e}\bar\psi^{1} dx\, \pi(z)dz \leq \left(\int \left(\int V^{1}\psi^{\e}\bar\psi^{1} dx\right)^2 \pi(z)dz\right)^{1/2} \leq ||V^{1}\psi^{\e}||_{\Gamma}\, 
||\psi^{1}||_{\Gamma}\,.
\end{align*}
Therefore, 
\[
    \frac{d}{dt} \| \psi^{1}\|_{\Gamma} \le\frac{1}{\e} \| V^1 \psi^{\e}\|^2_{\Gamma}\,. 
\]
For $t=O(1)$, this pessimistic estimate implies
\[  \| \psi^{1}\|_{\Gamma} =O\bigl(\e^{-1}\bigr). 
\]

{\color{black}{To summarize, in this part we emphasize the oscillatory behavior of the solution $\psi$ in the random space, which brings numerical challenges for the forward UQ problem. If one directly adopts the generalized polynomial chaos (gPC)-based Galerkin methods or stochastic collocation methods \cite{Xiu} to the semi-classical Schr\"odinger equation with random parameters, $\e$-dependent basis functions or quadrature points are needed to get an accurate approximation. There has been some work developed for this forward problem \cite{Liu20, GWPT}, where in our inverse problem case shares the similar difficulty. In the future work, to more efficiently sample from the random space, we will adopt numerical solvers that can resolve the $\e$-oscillations in the random variable. For simplicity of notations, we will omit the superscript $\e$ in $\psi^{\e}$ and use $\psi$ in the rest of the paper. }}

\section{Optimal control using neural networks} \label{sec:method}

\subsection{The time-splitting spectral method}\label{subsec:forward}

In the semiclassical regime where $\e \ll 1$, the solution to the Schr\"odinger equation \eqref{Schro-eqn}
is oscillatory both temporally and spatially, with an oscillation frequency of $O(1/\veps)$. This poses
tremendous computational challenges since one needs to numerically resolve, both spatially and temporally,
the small wave length of $O(\veps)$. The time-splitting spectral (TSSP) method, studied by Bao, Jin and Markowich in \cite{TSS}, is one of the most popular and highly accurate methods for such problems, where the meshing strategy $\Delta t=O(\e)$ and $\Delta x=O(\e)$ is required for moderate values of $\e$. Moreover, in order to just compute accurately the physical observables (such as position density, flux, and energy), one still needs to resolve the spatial oscillations, but the time step $\Delta t=o(1)$ is much more relaxed \cite{TSS, GJP, JMS}. Recently a rigorous uniform in $\e$ error estimate was obtained in \cite{Golse}, by using errors measured by a pseudo-metric in analogy to the Wasserstain distance between a quantum density operator and a classical density in phase space, with the regularity requirement for $V$ being 
$V \in C^{1,1}$. 

In this section, we review the first-order time-splitting spectral method studied in \cite[Section 2]{TSS}. 
Consider an one-dimensional spatial variable and a given potential $V(x)$. 
We choose the spatial mesh size $h=(b-a)/M$ for an even integer $M$, and the time step $k=\Delta t$, let the grid points and time step be
$$ x_j := a + j h, \qquad t_n:= n k, \qquad j=0,1, \cdots, M, \quad n=0,1,2, \cdots. $$
For the time discretization, from $t=t_n$ to $t=t_{n+1}$, the Schr\"odinger equation \eqref{Schro-eqn} is solved in the following two steps. First, one solves 
\begin{equation}\label{SP1} \e \psi_t - i \frac{\e^2}{2}\psi_{xx} = 0, \end{equation}
then \begin{equation}\label{SP2} \e \psi_t + i V(x) \psi = 0, \end{equation}
in the second step. 
We discretize \eqref{SP1} in space by the spectral method, then integrate in time {\it exactly}. Note that the ODE \eqref{SP2} can be solved exactly. 

Denote $\Psi_j^n$ by the numerical approximation of the analytic solution $\psi(t_n, x_j)$ to the Schr\"odinger equation \eqref{Schro-eqn}. 
Then the discretized scheme is given by
\begin{equation}
\label{SP-Scheme}
\begin{split}
&\Psi_j^{\ast} = \frac{1}{M} \sum_{l=-M/2}^{M/2-1} e^{-i \e k \mu_l^2/2}\, \hat\Psi_l^n \, e^{i \mu_l (x_j - a)}, \qquad
j = 0, 1, 2, \cdots, M-1, \\[6pt]
&\Psi_j^{n+1} = e^{-i V(x_j)k/\e} \Psi_j^{\ast}, 
\end{split}
\end{equation}
where the Fourier coefficients of $\Psi^n$ is defined as
$$ \hat \Psi_j^n = \sum_{j=0}^{M-1} \Psi_j^n\, e^{-i \mu_l (x_j - a)}, \qquad \mu_l = \frac{2\pi l}{b-a}, 
\quad l=-\frac{M}{2}, \cdots, \frac{M}{2}-1, $$
with $$ \Psi_j^0 = \psi(0, x_j), \quad j=0,1,2, \cdots, M. $$

We remark that instead of directly simulating the semi-classical Schr\"odinger equation, there are quite a few other methods which are valid in the limit $\varepsilon\to 0$, see \cite{BJM-Review} for a general discussion. In particular, many wave packets based methods have been introduced in past few years, which reduce the full quantum dynamics to Gaussian wave packets dynamics \cite{Heller}. In this work, we simply adopt the TSSP method as our deterministic solver in the learning algorithm

\subsection{Learning method for the control problem}\label{subsec:net}

Thanks to the nonlinear structure of deep neural network, it has shown great potential in approximating high dimensional functions and overcoming the curse of dimensionality. In recent years, deep learning has gained great success in solving high-dimensional PDEs, in both forward and inverse problem settings \cite{PINN, deepritz}. There have been studies that suggested learning-based methods on solving general control problems, such as \cite{Avila, Watter}. Recently, in \cite{sympocnet} the authors proposed SympOCnet to solve high dimensional optimal control problems with state constraints. The idea is to apply the Symplectic network, which can approximate arbitrary symplectic transformations, to perform a change of variables in the phase space and solve the forward Hamiltonian equation in the new coordinate system. In our work, we consider the control problem for the semiclassical Schr\"odinger equation and adopt neural networks to approximate the control field $V$ that may contain uncertainties. The neural network parameterized potential function is learnt by minimizing the discrepancies between the state solution of the system with neural network and the observation of the target state.   

In this section, we will describe the neural network structures under two different problem settings: (i) the deterministic case where the underlying target potential is fixed; (ii) the stochastic case where the target potential is parameterized by some random variables. In both problems, we will validate the efficiency of our proposed method by using both clean and noisy training data.

\subsubsection{Deterministic problem} 
\label{secsec:det}

\begin{figure}[!ht]
	\centering
	\includegraphics[width=0.9\textwidth]{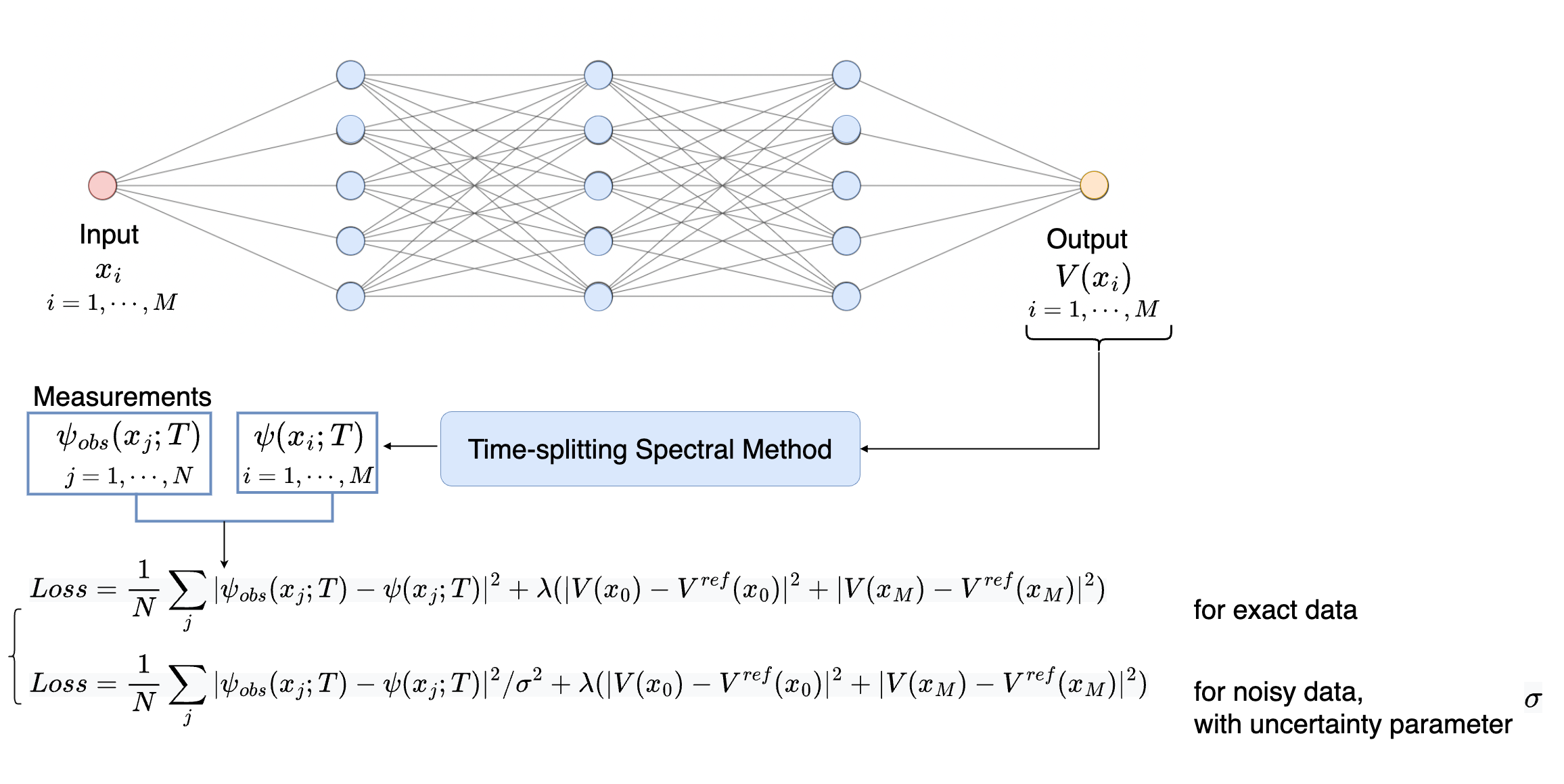}
 \vspace{-1cm}\caption{Illustration of the network for the deterministic problem.}\label{fig:det}
\end{figure}

In the deterministic problem, our goal is to learn a single target function $V(x)$ using the neural network. In this case, the input of the neural network is the spatial variable $\{x_k\}$, while the output is the value of the potential function at $x_k$, i.e., $\{V(x_k)\}$, $k=1,\cdots, M$. We will use $5$ fully connected layers with $50$ neurons per layer to build up the network. For the data points, assume the spatial domain $\Omega \in \mathbb{R}$ and temporal domain $[0,T]$, $N$ equally distributed points in $\Omega$ (where $N \ll M$) are taken, and the measurement data are the corresponding numerical solutions of the wave function at time $T$. This implies that the data pairs are chosen as $(x_i, \psi_{\text{obs}}(x_i))$ for $i=1, \cdots, N$ and $\psi_{\text{obs}}(x_i) \sim \mathcal{N}(\psi(x_i), \sigma^2)$. In our numerical examples, we set $N=50$ and $M=1000$. An illustration of the network for the deterministic problem is presented in Figure \ref{fig:det}. {\color{black} As noticed from Figure \ref{fig:det}, the input-output pairs for the fully connected neural network are $(x_i, V(x_i))$. The output of the neural network, i.e. the potential function, is then used to solve the forward Schr\"odinger equation by adopting the time-splitting spectral method. The predicted solution obtained at the final time step $\psi(x;T)$ is then compared with the measurement data $\psi_{\text{obs}}(x;T)$. The mismatch between the predicted solution and the measurement data will form the loss function. A pseudocode is presented in Algorithm \ref{alg:det}.}

\begin{algorithm} [!htb]
	\caption{Deterministic case}\label{alg:det}
	\begin{algorithmic}[1] 
		\Input{Neural network input $\{x_i\}_{i=1}^M$. Observation data $\{\psi_{\text{obs}}(x_j,T)\}_{j=1}^N$. Initialization of neural network parameters $\boldsymbol \theta_0$. }

		\For{For $k \gets 0: \#iterations$} 
		\State Get the output of the neural network $\{V(x_i;\boldsymbol \theta_k)\}_{i=1}^M$.
        \State Given $V(x;\boldsymbol \theta_k)$, solve equation (1.1) by time-splitting spectral method and get the solution $\psi(x,T;\boldsymbol \theta_k)$. 
        \State Compute the mismatch between $\psi_{\text{obs}}(x, T)$ and $\psi(x, T;\boldsymbol \theta_k)$, and get the loss.
        \State Use SGD type or SGLD method to update the network parameter and get $\boldsymbol \theta_{k+1}$.  
        \EndFor
     \Output{The solution of (1.1) $\psi(x_j, t_m)$ at all spatial locations and all time steps of interest.}
	\end{algorithmic}
\end{algorithm}

\subsubsection{Stochastic problem}
\label{secsec:sto}

In the stochastic problem, our goal is to learn a set of functions described by a stochastic potential function $V(x;z)$ containing a random parameter $z$, by training the DNN. We will utilize the DeepONet architecture developed in \cite{LJPZG21}. 

First we give a brief overview of DeepONet, which is a powerful tool designed to learn continuous nonlinear operators. Denote $G$ by an operator with input function $u$; for any coordinate $y$ in the domain of $G(u)$, the output $G(u)(y)$ is a number. DeepONet aims to approximate $G$ with a neural network $G_{\boldsymbol{\theta}}$ parameterized by $\boldsymbol{\theta}$, which takes inputs $(u,y)$ and returns the output $G(u)(y)$. The architecture of DeepONet is composed of a branch net and a trunk net. In the unstacked setting, the branch net encodes the discrete input function $u$ into the features represented by $[b_1, \cdots, b_q]$, and the trunk net takes the coordinate $y$ as input and encodes it into the features represented by $[t_1, \cdots, t_q]$. Then the dot product of $\textbf{b}$ and $\textbf{t}$ provides the final output of DeepONet, i.e.
 \begin{equation*}
	G_{\boldsymbol{\theta}}(u)(y) = \sum_{k=1}^q b_k(u(x_1), \cdots, u(x_N)) t_k(y). 
 \end{equation*} 
 The parameter $\boldsymbol{\theta}$ consists of all weights and biases in the branch and trunk net. 
 
\begin{figure}[!ht]
	\centering
	\includegraphics[width=0.8\textwidth]{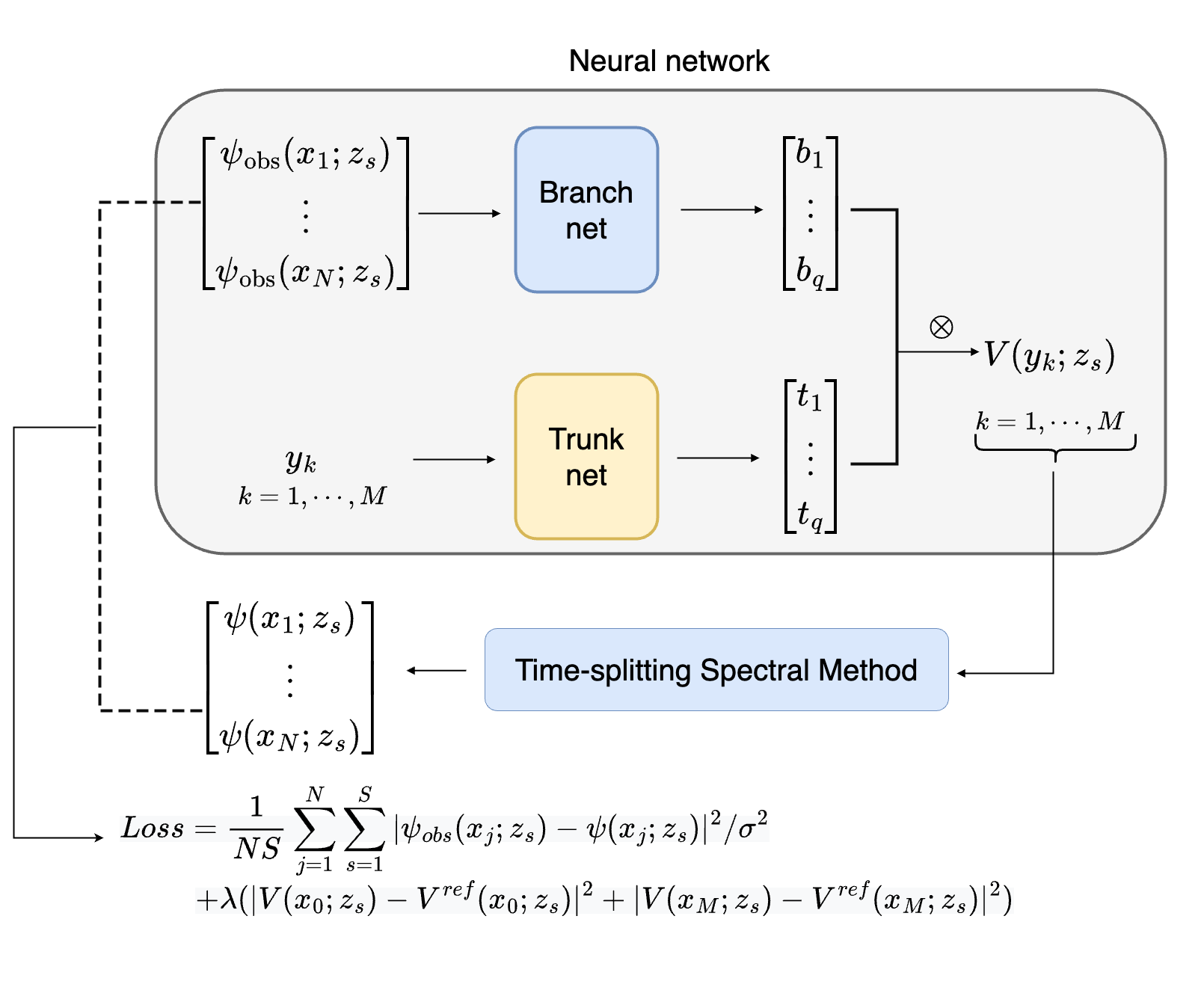} \vspace{-1cm}\caption{Illustration of the network for the stochastic problem.}\label{fig:sto}
\end{figure}

In our setting, we aim to approximate the parameterized potentials $V(x;z)$ using $G_{\boldsymbol{\theta}}$ that takes the discrete data $[\psi_{\text{obs}}(x_1; z),\cdots,\psi_{\text{obs}}(x_N; z)]$ and the coordinate $y_k$ as inputs. Here $k=1, \cdots, M$. We note that for each $z$, there are $N$ sensors that provide the observation data $\psi_{\text{obs}}(\cdot, z)$, thus the dataset size is equal to the product of $M$ and the number of $z$ samples. The value of $G_{\boldsymbol{\theta}}(\psi_{\text{obs}}(\cdot; z))(y_k)$ is a prediction of $V(y_k;z)$. Utilizing the predictions from the DeepONet, namely $V(y_k;z)$ ($k=1, \cdots, M$), the time-splitting spectral method is then applied to compute the value of wave functions $\psi(y_k, z)$. We aim to minimize the mismatch between the observations $\psi_{\text{obs}}(x_j, z)$ and the numerical solutions $\psi(x_j, z)$ at all sensor locations $x_j$ for all $z$. An illustration of the network for the stochastic problem is presented in Figure \ref{fig:sto}.  A pseudocode is presented in Algorithm \ref{alg:sto}.

\begin{algorithm} [!htb]
	\caption{Stochastic case}\label{alg:sto}
	\begin{algorithmic}[1] 
		\Input{Neural network input $\{\psi_{\text{obs}}(x_j, t_m; z_s)\}_{j=1}^N$ for some stochastic samples $z_s$, at few time instances $t_m$, as well as spatial points $\{y_k\}_{k=1}^M$. Observation data $\{\psi_{\text{obs}}(x_j,T; z_s)\}_{j=1}^N$. Initialization of neural network parameters $\boldsymbol \theta_0$. }

		\For{For $k \gets 0: \#iterations$} 
		\State Get the output of the neural network $\{V(y_k;z_s;\boldsymbol \theta_k)\}_{k=1}^M$.
        \State For each $V(x;z_s;\boldsymbol \theta_k)$, solve equation (1.1) by time-splitting spectral method and get the solutions $\psi(x, t;\boldsymbol \theta_k)$ at all spatial points and time instances. 
        \State Compute the mismatch between $\{\psi_{\text{obs}}(x_j,T; z_s)\}_{j=1}^N$ and $\psi(x_j,T;z_s;\boldsymbol \theta_k)$ (at the observational spatial and temporal points) over all samples of $z_s$, and get the loss.
        \State Use SGLD method to update the network parameter and get $\boldsymbol \theta_{k+1}$.  
        \EndFor
     \Output{For each $z_s$, the solution of (1.1) $\psi(x_j, t_m; z_s)$ at all spatial locations and all time steps of interest.}
	\end{algorithmic}
\end{algorithm}

\subsubsection{Training of the neural network} 
\label{secsec:sgld}

When dealing with large-scale problems, traditional Bayesian inference methods, e.g., Markov chain Monte Carlo (MCMC)\cite{MCMC} have shown disadvantages due to  extremely expensive computational cost of handling the whole dataset at each iteration. To tackle problems with large datasets, deep learning algorithms such as stochastic gradient descent (SGD) \cite{RM51} are favorable and have been popularly used, since one only needs to employ a small subset of samples randomly selected from the whole dataset at each iteration. To bring together advantages of these two types of methods, Welling and Teh \cite{sgld} first proposed the stochastic gradient Langevin dynamics (SGLD) (also known as stochastic gradient MCMC) method. It adds a suitable amount of noise to the standard SGD and uses mini-batches to approximate the gradient of loss function. With the help of decreasing training step size $\eta_k$, it has demonstrated powerful and provided a transition between optimization and Bayesian posterior sampling \cite{sg-mcmc-convergence}.

We now briefly review the SGLD method. Denote $D = \{d_i\}_{i=1}^N =\{(\textbf{x}_i, \textbf{y}_i) \}_{i=1}^N $ by a given dataset, where ${\color{black}{\textbf{x}_i}}$ is the input and ${\color{black}{\textbf{y}_i}}$ is the corresponding noisy output. We let $\mathcal{NN}$ be a neural network parameterized by the parameter $\boldsymbol \theta$; the goal of its training is to find suitable parameters $\boldsymbol \theta$ such that $F(\mathcal{NN}(\textbf{x}_i; \boldsymbol \theta)) \approx \textbf{y}_i$ ($i=1,\cdots, N$). Due to the noise in measurement data, we assume the parameters are associated with uncertainties and obey a prior distribution $p(\boldsymbol \theta)$. The uncertainties in the parameters $\boldsymbol \theta$ can be captured through Bayesian inference to avoid overfitting. Let $d^j$ be a mini-batch of data with size $n$, the likelihood can be written as 
 \begin{equation*}
 	p (d^j| \boldsymbol \theta) = \frac{1}{(2\pi \sigma^2)^{n/2}} \exp \Big\{-\frac{  \sum \limits_{\textbf{x}^j_i \in d^j} (\textbf{y}^j_i  -F(\mathcal{NN}(\textbf{x}^j_i ; \boldsymbol \theta))  )^2 }{ 2 \sigma^2 } \Big\},
 \end{equation*}
where $\sigma$ is standard deviation of the Gaussian likelihood. In our case, for the dataset $d^j = (\textbf{x}^j_i, \textbf{y}^j_i)$, $\textbf{x}_j^i$ corresponds to the input $[\psi_{\text{obs}}(x_1; z),\cdots,\psi_{\text{obs}}(x_N; z), y]$, $\textbf{y}_j^i$ corresponds to the labels $[\psi_{\text{obs}}(x_1; z),\cdots,\psi_{\text{obs}}(x_N; z)]$ and $F$ maps the output of the neural network output $\mathcal{NN}(\textbf{x}_i; \boldsymbol \theta)$ which approximates $V(y,z)$ to the quantities of interest $\psi(y;z;T)$ with $T$ the final simulation time. According to the Bayes' theorem, the posterior distribution of $\boldsymbol \theta$, given the data $D$, then follows $p(\boldsymbol \theta|D) \propto p(\boldsymbol \theta) \prod_{i=1}^{N} p(d_i|\boldsymbol \theta)$. 
 
To sample from the posterior, one efficient proposal algorithm is to use the gradient of the target distribution. Let $\eta_k$ be the learning rate at epoch $k$ and $\tau> 0$ be the inverse temperature, the parameters will be updated by SGLD based on the following rule: 
\begin{equation*}
{\boldsymbol \theta}_{k+1} = {\boldsymbol \theta}_{k} + \eta_k \nabla_{\boldsymbol \theta} \tilde{L} ({\boldsymbol \theta}_{k} ) + \mathcal{N}(0, 2\eta_k \tau^{-1}). 
\end{equation*}
Here for a subset of $n$ data points $d^j = \{d_1^j, \cdots, d_n^j \}$, 
\begin{equation*}
\nabla_{\boldsymbol \theta} \tilde{L} (\boldsymbol \theta ) = \nabla_{\boldsymbol \theta} \log p(\boldsymbol \theta) + \frac{N}{n} \sum_{i=1}^n \nabla_{\boldsymbol \theta} \log p(d^j_i|\boldsymbol \theta)
\end{equation*} 
is the stochastic gradient computed by using a minibatch that approximate the true gradient of the loss function $\nabla_{\boldsymbol \theta} {L} (\boldsymbol {\theta})$.

However, if the components of the network parameters $\boldsymbol \theta$ have different scales, the invariant probability distribution for the Langevin equation is not isotropic.
If one still uses a uniform learning rate in each direction, this may leads to slow mixing \cite{chen2014stochastic, dauphin2015equilibrated, dauphin2014identifying, PSGLD, psgld-sa, zhang2011quasi}. To incorporate the geometric information of the target posterior, stochastic Gradient Riemann Langevin Dynamics (SGRLD) \cite{SGRLD} generalizes SGLD on a Riemannian manifold. Consider the probability model on a Riemann manifold with some metric tensor $P^{-1}(\boldsymbol \theta)$, in SGRLD, the parameter is updated at the $k$-th iteration by the following rule: 
\begin{equation} \label{eq:sgrld}
{\boldsymbol \theta}_{k+1} = {\boldsymbol \theta}_{k} + \eta_k  \left[ P(\boldsymbol \theta_k ) \nabla_{\boldsymbol \theta} \tilde{L} ({\boldsymbol \theta}_{k} ) + \Gamma({\boldsymbol \theta}_k)  \right] +  \mathcal{N}(0, 2\eta_k \tau^{-1}P(\boldsymbol \theta_k )  )
\end{equation}
where $\displaystyle\Gamma_i({\boldsymbol \theta}_k) = \sum_j \frac{\partial P_{ij}(\boldsymbol \theta_k)  }{\partial \theta_j}$. One popular and computationally efficient approach to approximate $P(\boldsymbol \theta_k )$ is to use a diagonal preconditioning matrix \cite{PSGLD, psgld-sa}, 
that is, 
\begin{align}\label{eq:preconditioner}
	P(\boldsymbol { \theta}_k) &=  diag^{-1} (\lambda + \sqrt{V(\boldsymbol { \theta}_k)}),  \\
	V(\boldsymbol { \theta}_k) &=  (1-\omega_k) V(\boldsymbol { \theta}_{k-1}) + \omega_k g(\boldsymbol { \theta}_k) \circ g(\boldsymbol { \theta}_k), 
\end{align}
where $\lambda$ is a regularization constant, $g(\boldsymbol {\theta}_k) = \nabla_{\boldsymbol {\theta}} \tilde{L}(\boldsymbol {\theta}_k)$ is the stochastic gradient, the operator $\circ$ denotes a elementwise multiplication, and $\omega_k \in (0,1)$ is a weight parameter used in the moving average $V(\boldsymbol { \theta}_k)$.
In our framework, we will use the preconditioned SGLD to train the network parameters.

\section{Numerical results}\label{sec:numerical}

In our numerical experiments, we consider two types of potential functions, the deterministic and stochastic potential. In the deterministic case, the potential $V$ is only spatially dependent. In the stochastic problem, the potential function $V(\cdot, z)$ is assumed to depend on a random parameter characterized by $z$. In particular, we consider a simple example with 
$V(x,z) = (1+0.5z) x^2$, where $z$ is a random variable following the uniform distribution in $[-1,1]$. 

\subsection{Test I: A Deterministic Potential}

In the first problem setup, we assume the potential function as $V(x) = x^2$. The network architecture introduced in Section \ref{secsec:det} is adopted, and we train the network by using standard SGD and SGLD studied in Section \ref{secsec:sgld}. 
For the observation data, we choose it to be the electron wave function $\psi$ solved by the Schr\"odinger equation \eqref{Schro-eqn} at several spatial locations and time instances using the forward TSSP solver, given the reference potential function $V$.

We first consider that there is no noise in the observation data and apply both SGD and SGLD to train. The numerical results show that the wave function $\psi$ obtained from the network of both training algorithms matches well with the observation data, while it is also noticeable that the SGLD gives a slightly better approximation of the potential function. We then consider when some noise is added to the observation data, one can just apply SGLD to train the network in order to more accurately capture the uncertainties in the target potential function. 

	\subsubsection{$V(x) = x^2$, no noise in the observation and by SGD}
	
\begin{figure} 
		\centering
		\begin{subfigure}{0.45\textwidth}
			\centering
			\includegraphics[width=\textwidth]{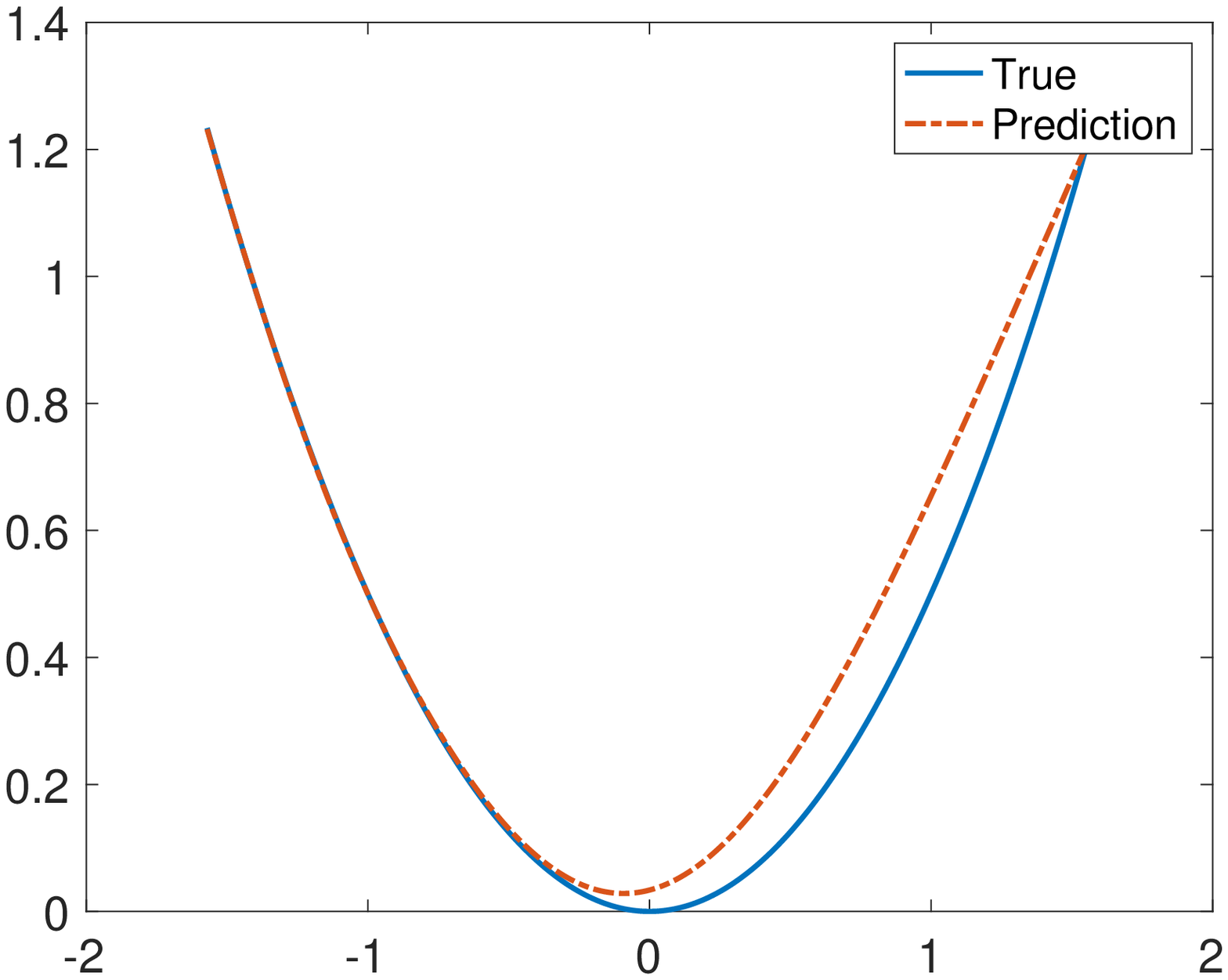}
			\caption{}
		\end{subfigure}
		\begin{subfigure}{0.45\textwidth}
			\centering
			\includegraphics[width=\textwidth]{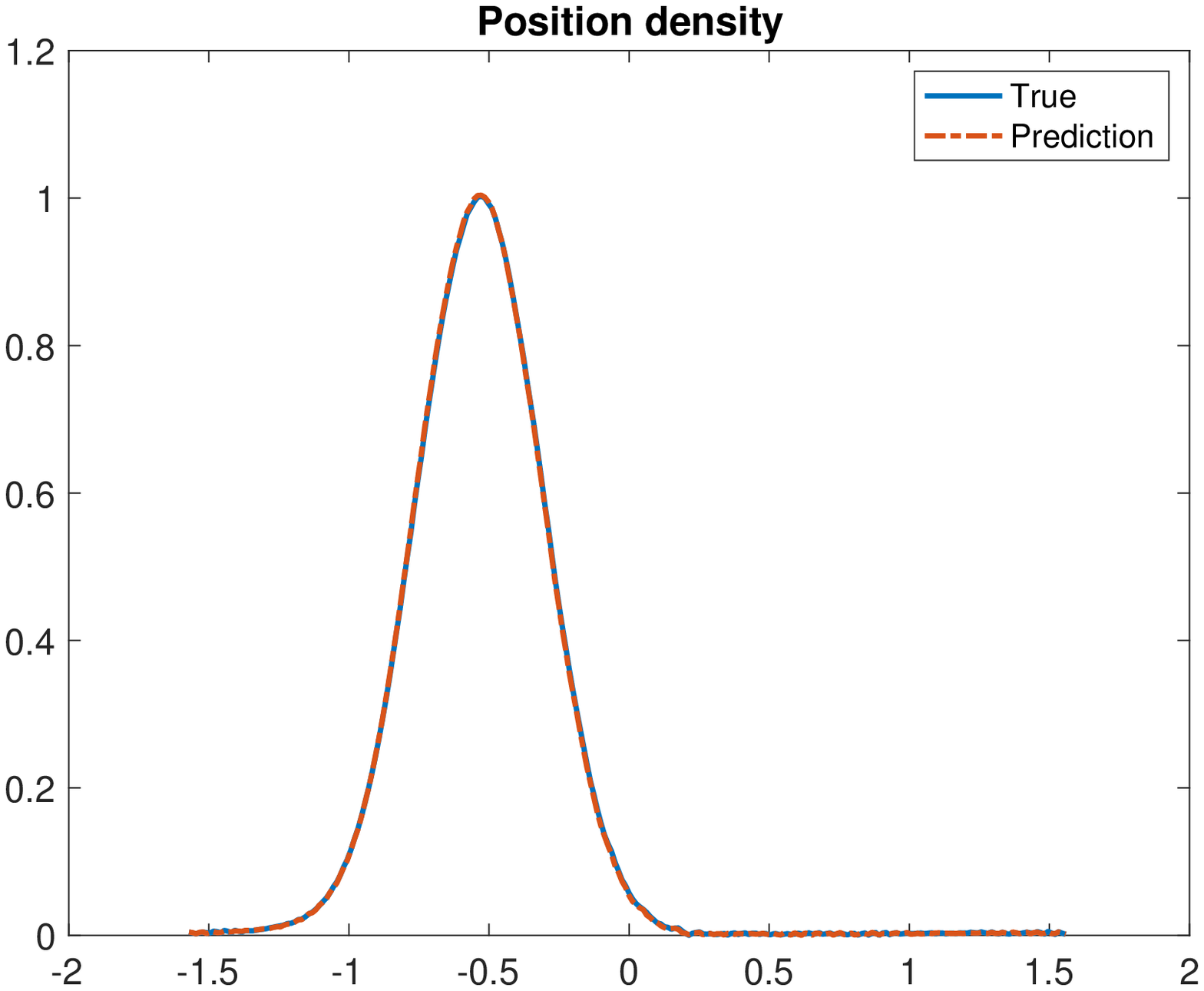}
			\caption{}
		\end{subfigure}
        \begin{subfigure}{\textwidth}
          \centering
          \includegraphics[width=1.0\textwidth]{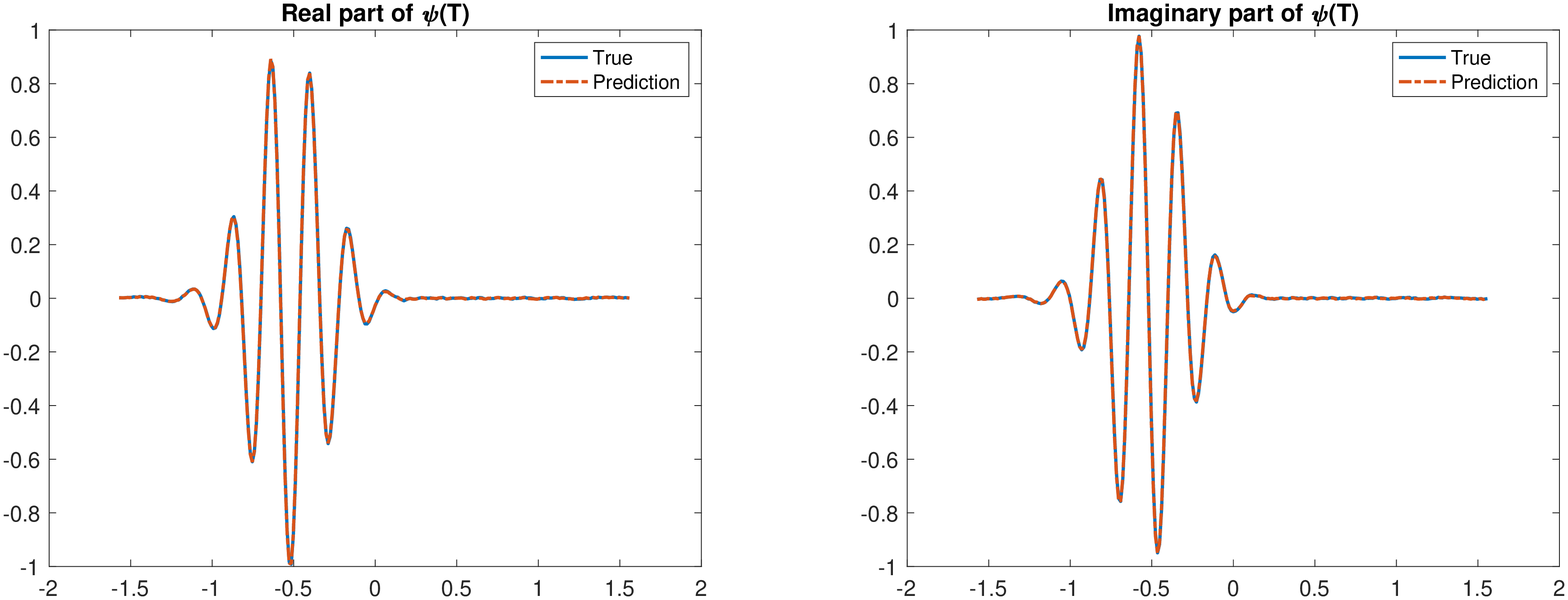}
          \caption{}
		\end{subfigure}
		\caption{Test I case 1: $V(x) = x^2$, without noise in the observation and by SGD. (a) True and predicted value of the potential function. (b) True and predicted value of the position density at time $T$. (c) True and predicted value of the wave function at time $T$. }\label{fig:ex1_sgd}
	\end{figure}

In this case, we let the reference potential function be $V(x)=x^2$, here the observation data is clean and without noise interference, SGD method is used in our training algorithm. {\color{black} In the forward solver, the spatial mesh size is $\pi/250$ and the temporal mesh size is $6.25\times 10^{-4}$. The learning rate is $10^{-4}$ and the total training epoch is $20000$.} In Figure \ref{fig:ex1_sgd} (a), a comparison between the reference and predicted potential function obtained from the neural network is shown. We observe that there is some mismatch in the region when $x>0$, while the underlying reason remains to be discovered. 
 In Figure \ref{fig:ex1_sgd} (b)-(c), a comparison between the reference with predicted position density $n^{\varepsilon}$ and the wave function $\psi^{\varepsilon}$ (real and imaginary parts) is presented. We conclude that the predicted wave and density functions at the final time $T$, which are computed by solving the Schr\"odinger equation \eqref{Schro-eqn} under the neural network's predicted output potential, can provide good approximations to the solution quantities obtained by using the true potential $V(x)=x^2$ in the TSSP solver. 
	
\subsubsection{$V(x) = x^2$, no noise in the observation and by SGLD}

\begin{figure} 
		\centering
		\begin{subfigure}{0.45\textwidth}
			\centering
			\includegraphics[width=\textwidth]{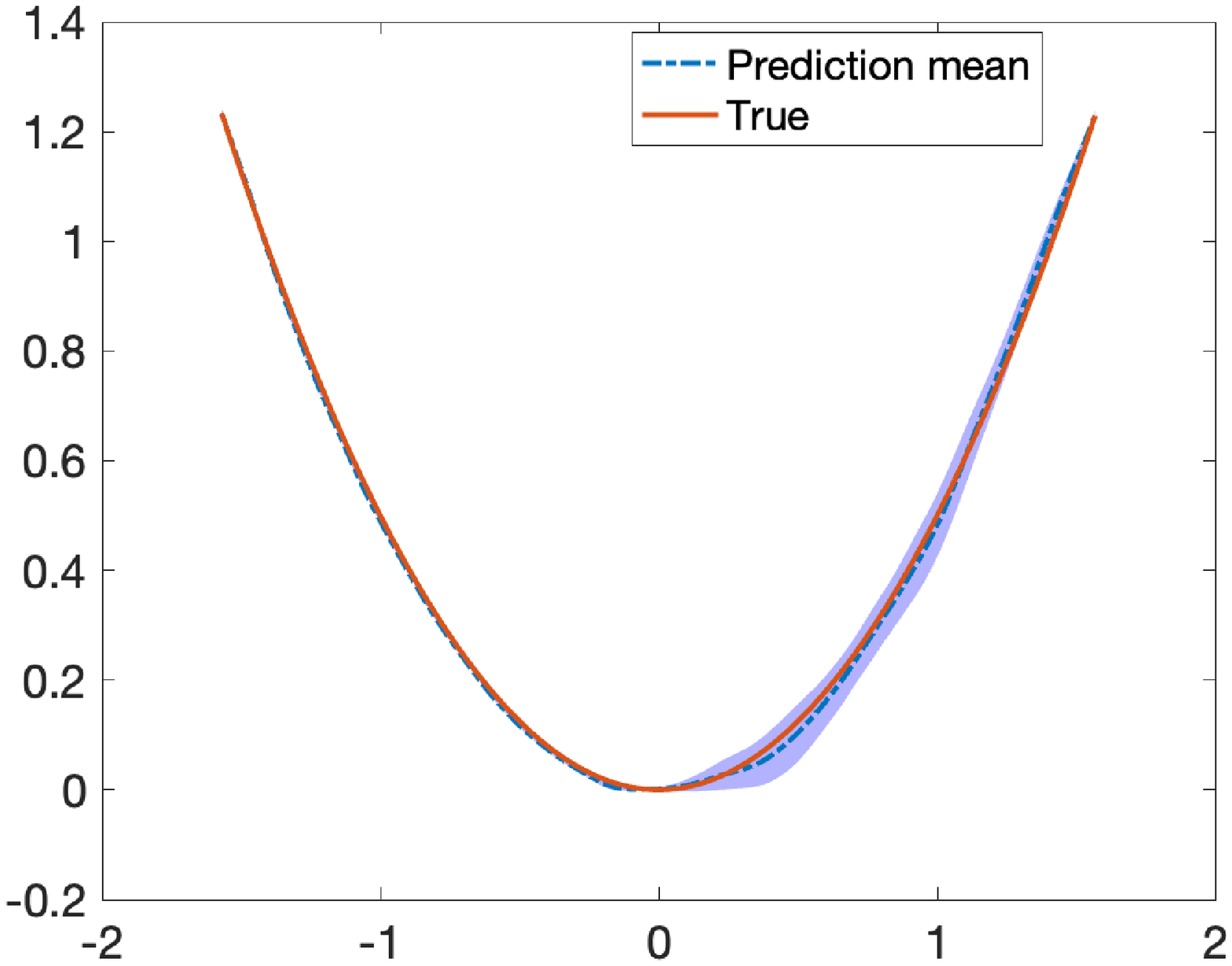}
			\caption{}
		\end{subfigure}
		\begin{subfigure}{0.45\textwidth}
			\centering
			\includegraphics[width=\textwidth]{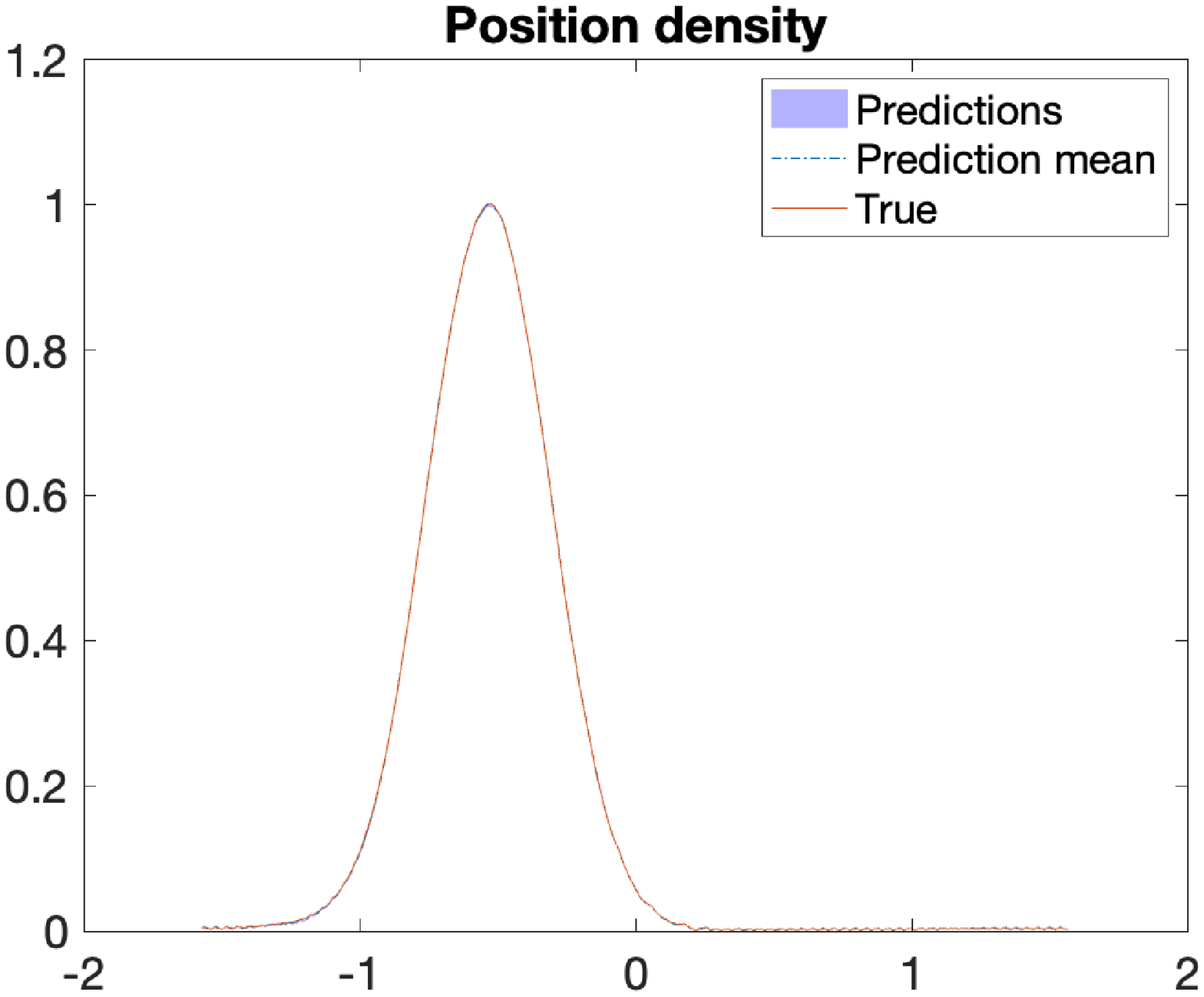}
			\caption{}
		\end{subfigure}
        \begin{subfigure}{\textwidth}
          \centering
          \includegraphics[width=1.0\textwidth]{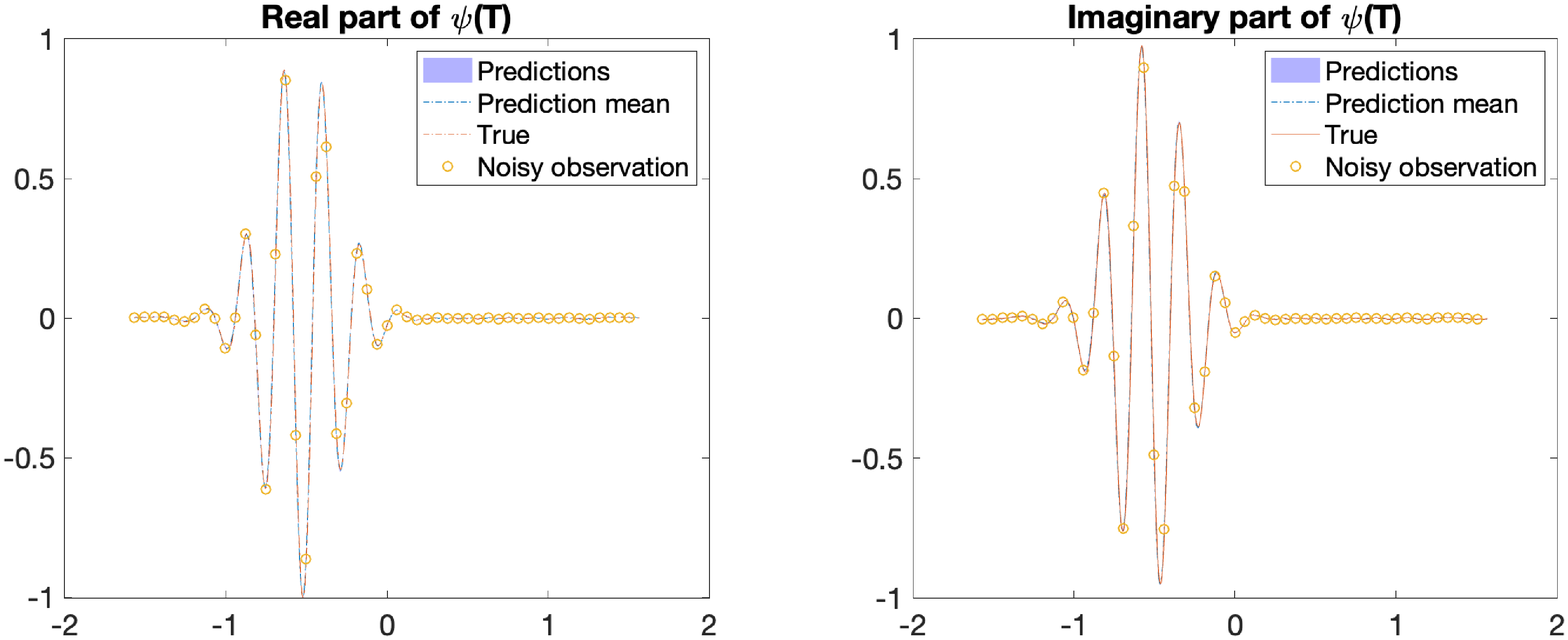}
          \caption{}
		\end{subfigure}
		\caption{Test I case 2: $V(x) = x^2$, without noise in the observation and by SGLD. (a) True and predictions of the potential function. (b) True and predictions of the position density at time $T$. (c) True and predictions the wave function at time $T$. }\label{fig:ex1_sgld}
	\end{figure}

In the second case, the problem setup is the same as the previous case, while we apply SGLD algorithm to train the neural network. {\color{black} In the forward solver, the spatial mesh size is $\pi/1000$ and the temporal mesh size is $3\times 10^{-3}$. The learning rate is $10^{-5}$ and the total training epoch is $10000$.} A comparison between the reference and predicted potential function is shown in Figure \ref{fig:ex1_sgld} (a). According to the nature of SGLD, we collect samples of neural network's parameters during the training process, then compute the mean and standard deviation of output potential functions (at each spatial point) obtained by using those parameter samples. The blue dashed line represents the mean of the predicted potential $V$, and the confidence interval are depicted by the shaded blue area in Figure \ref{fig:ex1_v_sgld}. Based on these two tests, we observe that SGLD provides more reliable results compared to the standard SGD, and the uncertainty is neglible in the prediction since the data is clean. 

In Figure \ref{fig:ex1_sgld} (b)-(c), we again present a comparison between the reference and predicted wave function $\psi^{\varepsilon}$ or position density $v^{\varepsilon}$ that is computed by the TSSP solver by using the predicted mean value of the potential. Similar to the previous test, it is obvious that the predicted wave or density can provide quite good approximations to the true data, i.e., the numerical solution at final time $T$ obtained by using the true potential $V(x)=x^2$ in the TSSP solver. 

\subsubsection{$V(x) = x^2$, noisy data and by SGLD}

\begin{figure} 
		\centering
		\begin{subfigure}{0.45\textwidth}
			\centering
			\includegraphics[width=\textwidth]{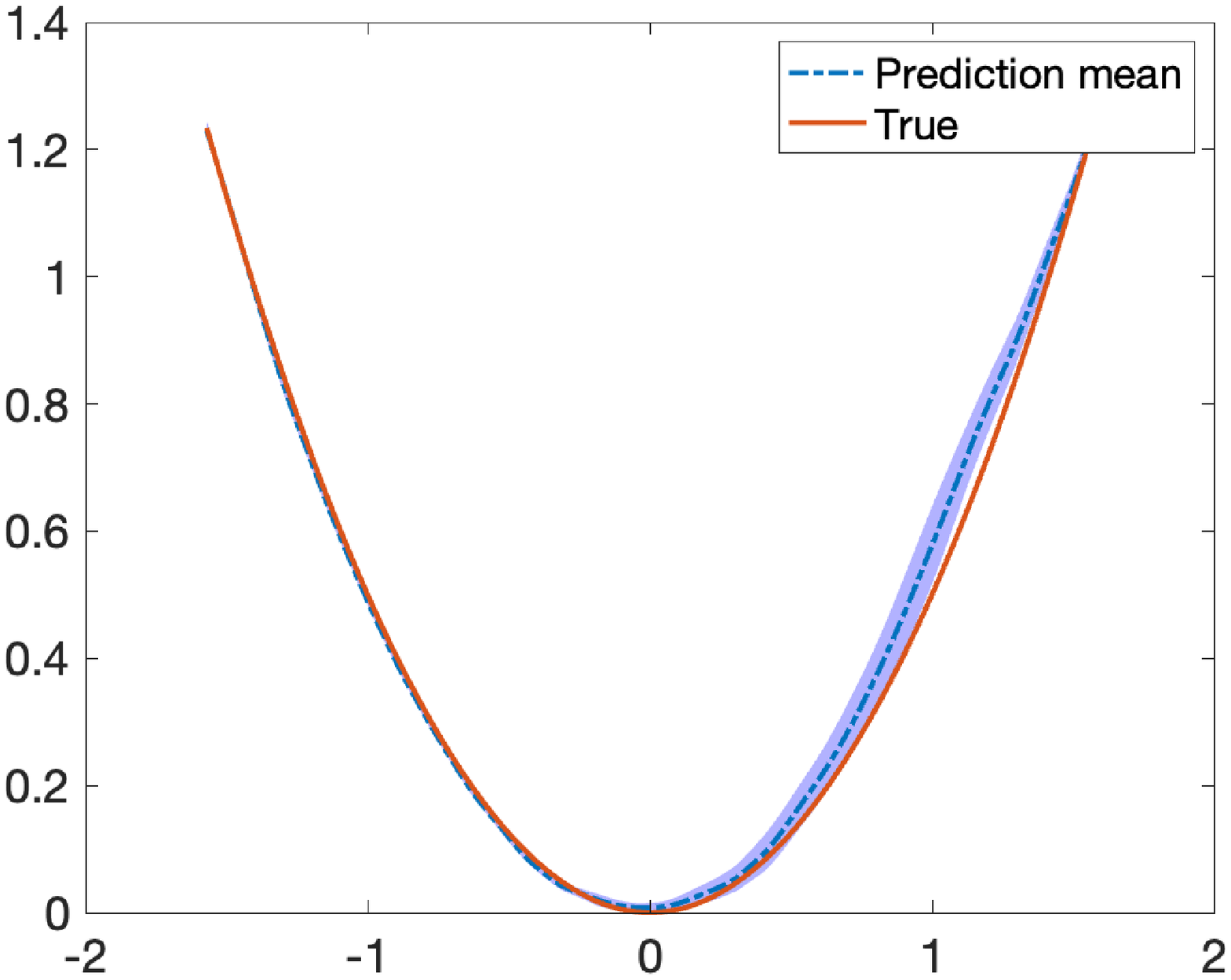}
			\caption{}
		\end{subfigure}
		\begin{subfigure}{0.45\textwidth}
			\centering
			\includegraphics[width=\textwidth]{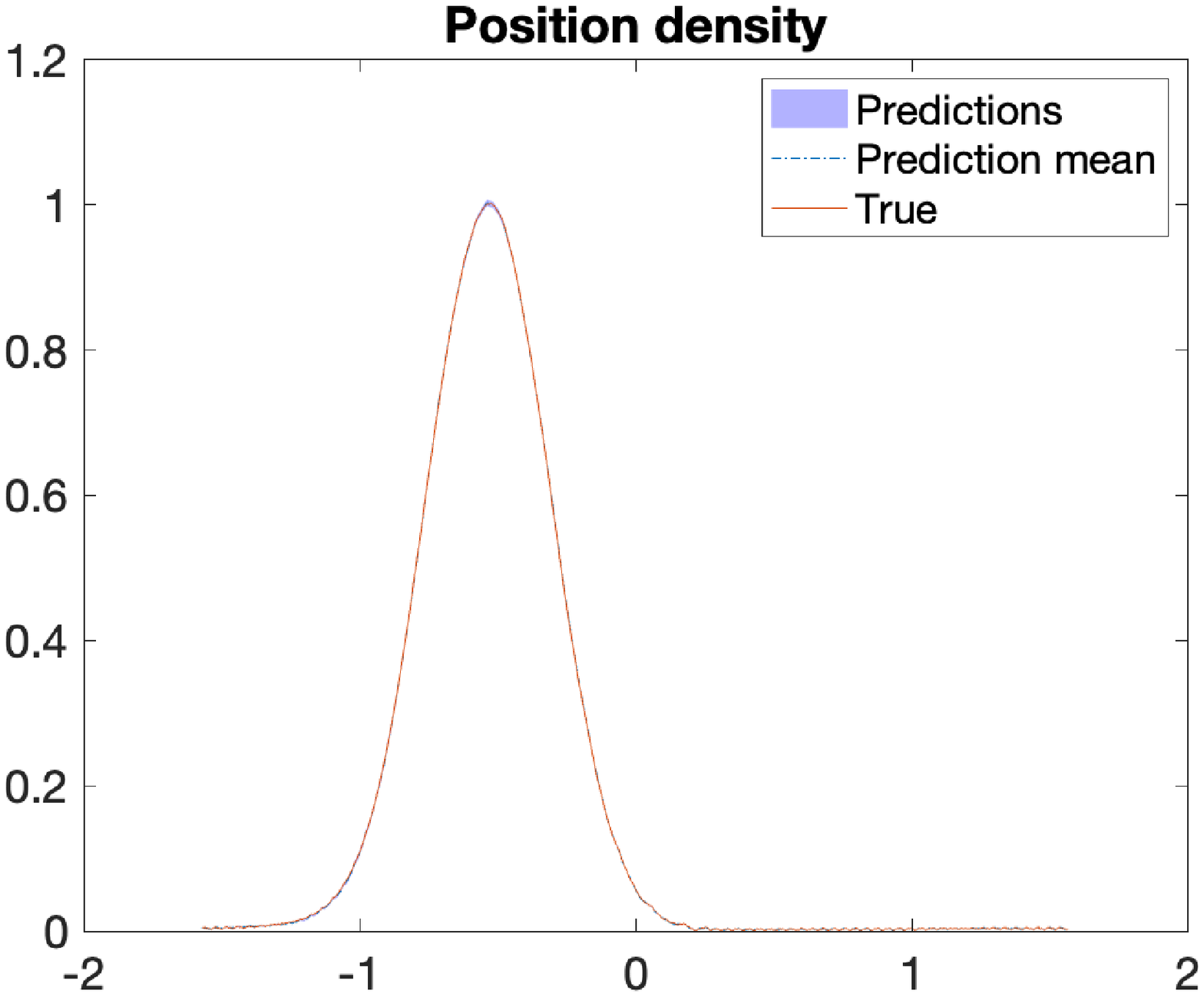}
			\caption{}
		\end{subfigure}
        \begin{subfigure}{\textwidth}
          \centering
          \includegraphics[width=1.0\textwidth]{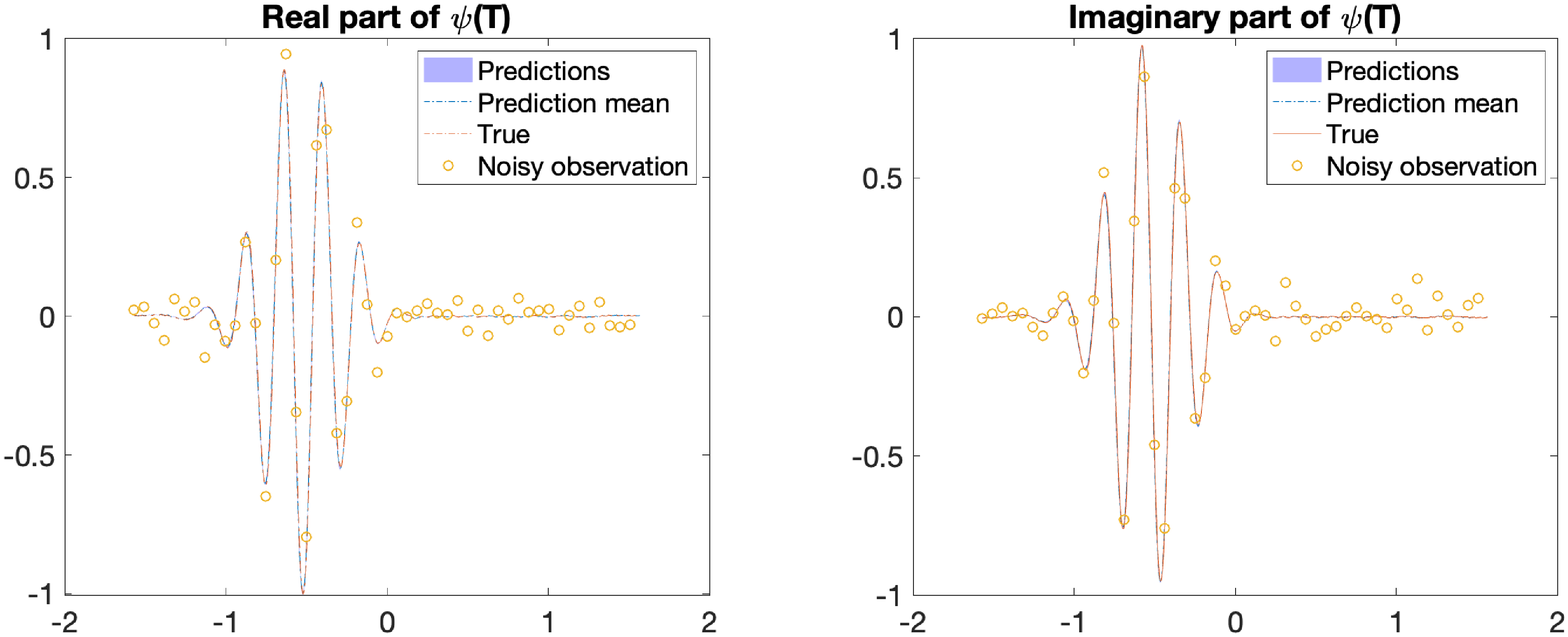}
          \caption{}
		\end{subfigure}
		\caption{Test I case 3: $V(x) = x^2$, noisy data and by SGLD. (a) True and predictions (with confidence interval) of the potential function. (b) True and predictions of the position density at time $T$. (c) True and predictions the wave function at time $T$. }\label{fig:ex1_sgld_n5}
	\end{figure}

In the thrid case, we consider some noise in the observation data and use SGLD to train the network. {\color{black}The mesh size in the forward solver, the learning rate and the training epochs are the same as in the previous subsection.}  We let the noise be a random variable that follows the normal distribution with mean $0$ and standard deviation $0.05$. In Figure \ref{fig:ex1_sgld_n5} (c), the yellow circles are the noisy values of $\psi^{\e}$ at $50$ equally spaced locations. A comparison between the reference, i.e., $V(x) = x^2$, with the predicted mean of the potential function is shown in Figure \ref{fig:ex1_sgld_n5} (a). One can observe that the predicted mean value is consistent with the reference potential, and the blue shaded area indicates that there are some uncertainties due to the noisy data, compared to the previous tests where there is no noise in the observation. 

Similarly, we can see from Figure \ref{fig:ex1_sgld_n5} (b)-(c), the predicted wave and density at final time $T$ that are computed using the mean of network's predicted potential $V$, capture well the true solution obtained by using $V(x)=x^2$ in the TSSP solver. Therefore, we conclude that SGLD can deal with the noisy data and provide reliable results. 

\subsection{Test II: A Stochastic Potential}

In Test II, we consider a stochastic potential, $V(x,z) = (1+0.5z) x^2$, where $z$ follows the uniform distribution on $[-1,1]$. 
To generate the dataset, we first take eight Gauss-Legendre points for $z\in [-1,1]$. For each $z_k$ ($k=1, \cdots, K$), i.e., each specific potential $V(x;z_k)$, we have the corresponding noisy measurement data $\psi_{\text{obs}}(x;z_k)$ at the final time instance $T=0.6$. The observation $\psi_{\text{obs}}(x;z_k) \sim \mathcal{N}(\psi(x;z_k), \sigma^2)$ where $\sigma$ is the standard deviation. The wave functions at the final time instance $\psi(x;z_k)$ is computed using the time-splitting spectral method on a $640\times 1000$ temporal-spatial grid. Then for each $z_k$ we select $N$ sensor locations to collect the measurement data, the sensors are uniformly located in the spatial domain $\Omega = [-\pi/2, \pi/2]$. We will take $N=20, 50$ in the numerical tests. {\color{black} In the forward solver, the spatial mesh size is $\pi/1000$ and the temporal mesh size is $6.25\times 10^{-4}$. The learning rate is $10^{-5}$ and the total training epoch is $10000$.}

The input of the network then consists of the spatial evaluation point $x_i$, and the real part $\Re(\psi_{\text{obs}}(x_1;z_k)), \cdots, \Re(\psi_{\text{obs}}(x_N;z_k)$ and the imaginary part $\Im(\psi_{\text{obs}}(x_1;z_k)), \cdots, \Im(\psi_{\text{obs}}(x_N;z_k)$ of the observation data. The output of the network is the value of potential at $x_i$, i.e., $V(x_i;z_k)$. The number of training samples is equal to the product of $M$ (the number of evaluation points $x_i$) and the number of $z$ samples. We assume that the values of $V(x,z)$ at the endpoints $x=-\frac{\pi}{2}$ and $x=\frac{\pi}{2}$ are known for the training samples. The loss function consists of three parts, (1) the mismatch between the observation data $\psi_{\text{obs}}(x;z_k)$ and the $\psi$ computed using neural network predicted potential function, (2) the mismatch between the true potential and neural network predicted potential at the endpoints of the spatial domain, and (3) a regularization term on the potential. After training, we will obtain the full potential profile for different $z$ samples. In the testing stage, we will only have noisy observations of the wave function at final time $T$ without knowing any information of the true potential function. We will feed the a set of spatial location $x_i$ as well as the observation data into the neural network, and obtain the predictions of the potential evaluated at these points $x_i$.
 
\begin{figure}[!ht]
	\centering
	\includegraphics[width=0.45\textwidth]{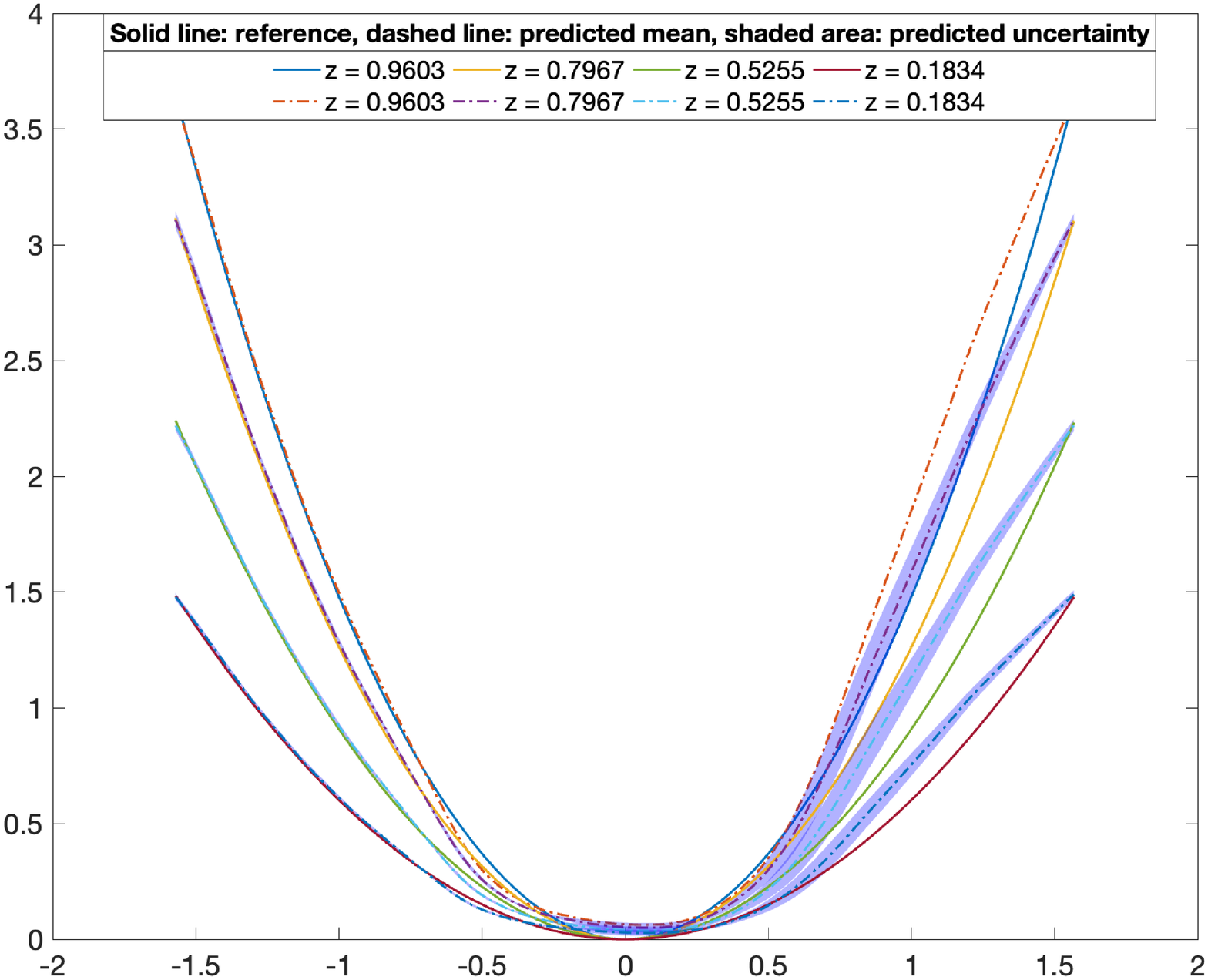}
 \includegraphics[width=0.45\textwidth]{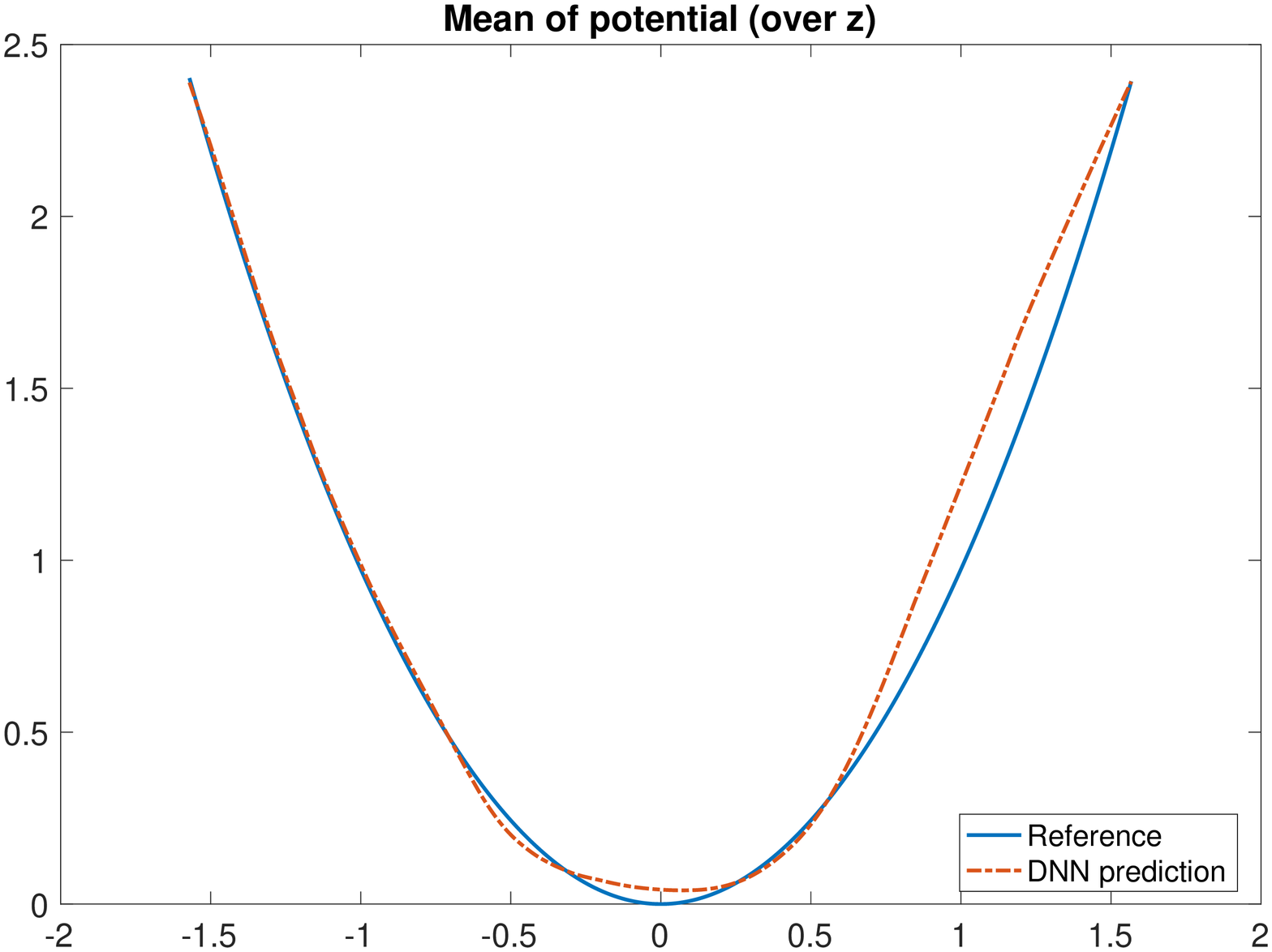}
	\caption{Test II, true and predicted value of the potential function $V(x;z) = (1+0.5z )x^2$ when the number of sensors is $50$. Left: different $z$s, right: mean prediction with respect to $z$.}
	\label{fig:sto_V_20}
\end{figure}

\begin{figure}[!ht]
	\centering
	\includegraphics[width=0.9\textwidth]{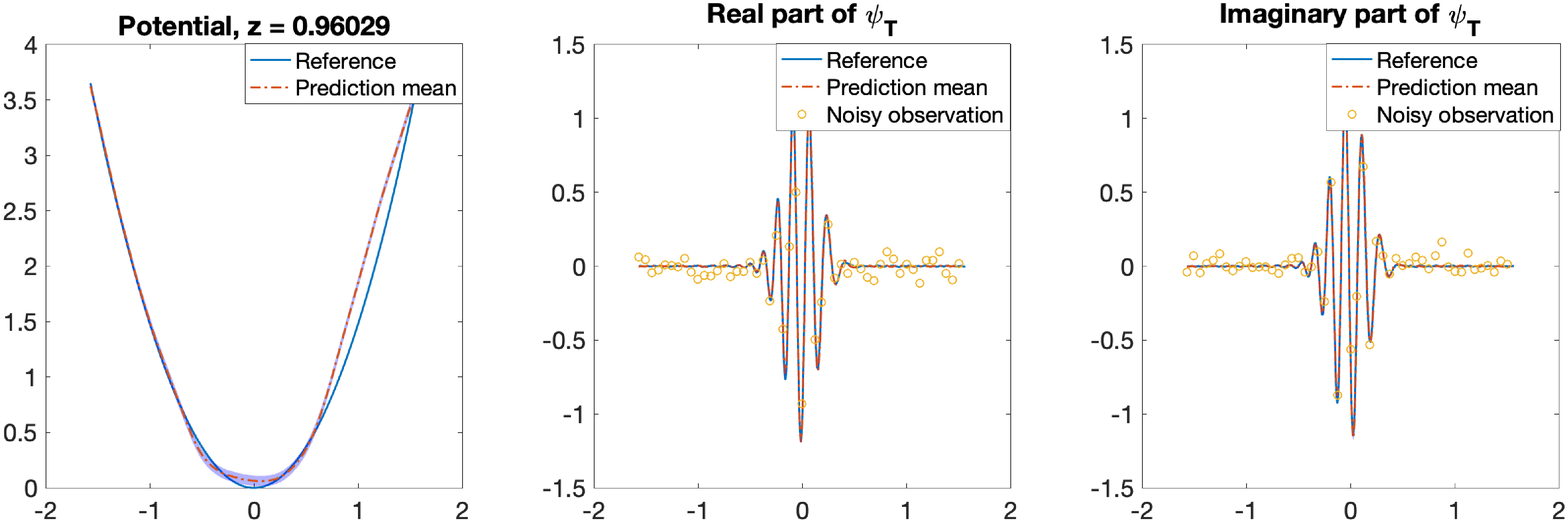}
	\caption{Test II, true and predicted value of the potential function $\psi$ at final time $T=0.6$, for a training sample $z=0.9603$, 50 sensors.}
	\label{fig:sto_psi1_20}
\end{figure}

\begin{figure}[!ht]
	\centering
	\includegraphics[width=0.9\textwidth]{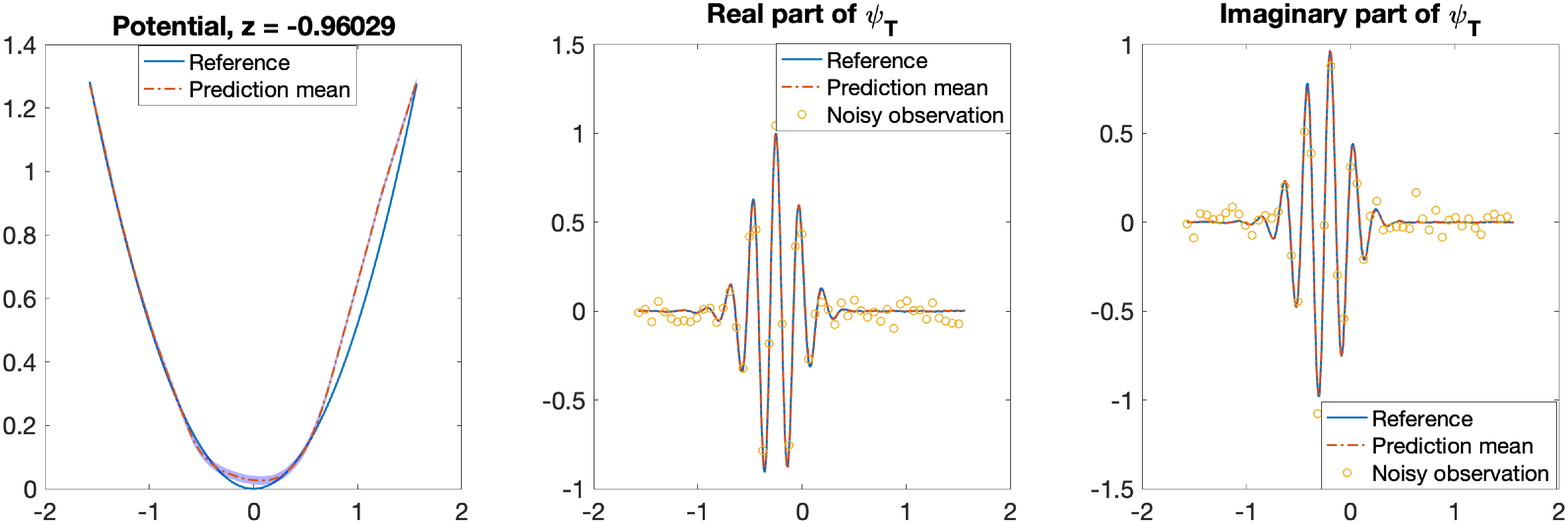}
	\caption{Test II, true and predicted value of the potential function $\psi$ at final time $T=0.6$, for a training sample $z=-0.9603$, 50 sensors.}
	\label{fig:sto_psi2_20}
\end{figure}

\begin{figure}[!ht]
	\centering
	\includegraphics[width=0.9\textwidth]{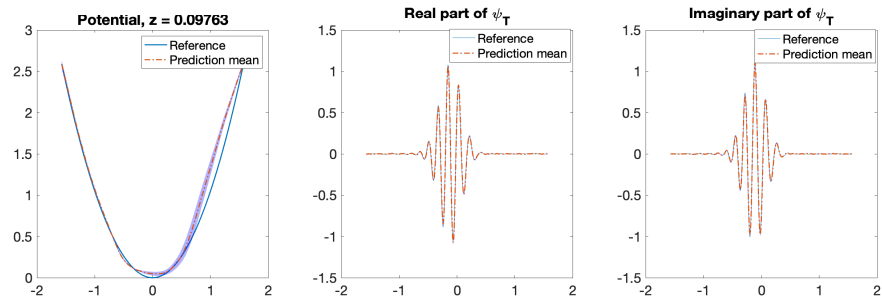}
	\caption{Test II, testing case: $z=0.0976$. True and predicted value of the potential function $\psi$ at final time $T=0.6$, 50 sensors.}
	\label{fig:sto_psi_test_20}
\end{figure}

We first show the predictions of $V(x;z)$ for some training samples of $z$ when there are $50$ sensors and $\psi_{\text{obs}}(x;z) \sim \mathcal{N}(\psi(x;z), 0.05)$. The comparison of predictions and references of $V(x;z) =(1+0.5z) x^2$  for four different $z$ values ($z = [0.9603, 0.7967, 0.5255, 0.1834]$) are presented on the left of Figure \ref{fig:sto_V_20}. The expected value of $V$ over the random variable $z$ are computed using 8 Legendre quadrature points in the interval $z\in [-1,1]$, and the comparison of the predicted mean and reference mean are shown on the right of Figure \ref{fig:sto_V_20}.
With large numbers of observation data and suitable amount of noise in the data, the neural network can provide reasonable approximations for the potential functions. The corresponding predictions of wave function $\psi$ (computed using predicted potential functions) at the final time $T=0.6$ with different values of $z$ are shown in Figures \ref{fig:sto_psi1_20}, \ref{fig:sto_psi2_20}. We observe good agreements between the predictions and the true values of the wave functions. A testing case for $z=0.0976$ is shown in Figure \ref{fig:sto_psi_test_20}. It shows that our trained neural network can generalize well to new samples of $z$.

\begin{figure}[!ht]
	\centering
	\includegraphics[width=0.45\textwidth]{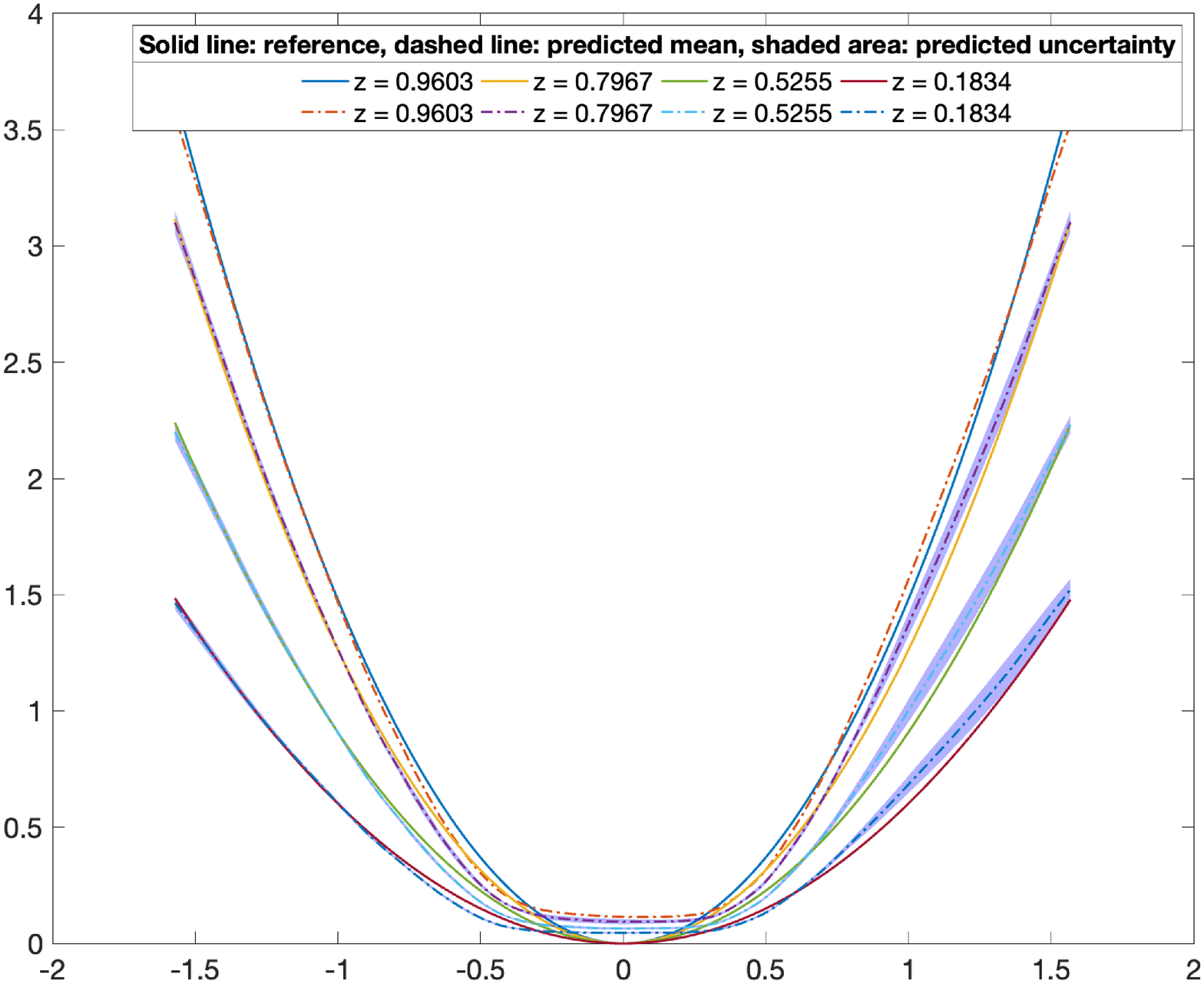}
 \includegraphics[width=0.45\textwidth]{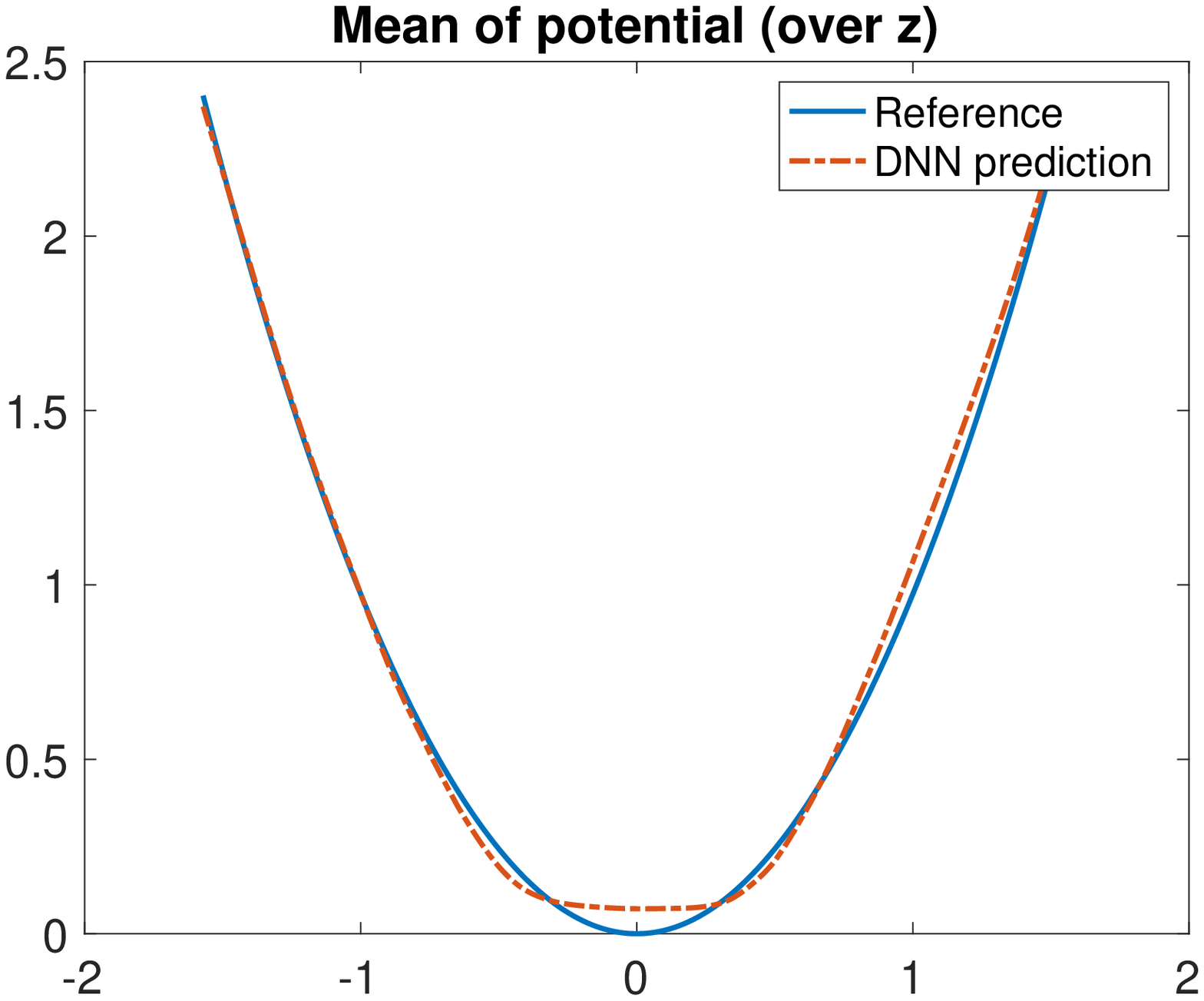}
	\caption{Test II, 20 sensors, true and predicted value of the potential function $V(x;z) = (1+0.5z )x^2$. Left: different $z$s, right: mean prediction with respect to $z$.}
	\label{fig:sto_V_50}
\end{figure}

\begin{figure}[!ht]
	\centering
	\includegraphics[width=0.9\textwidth]{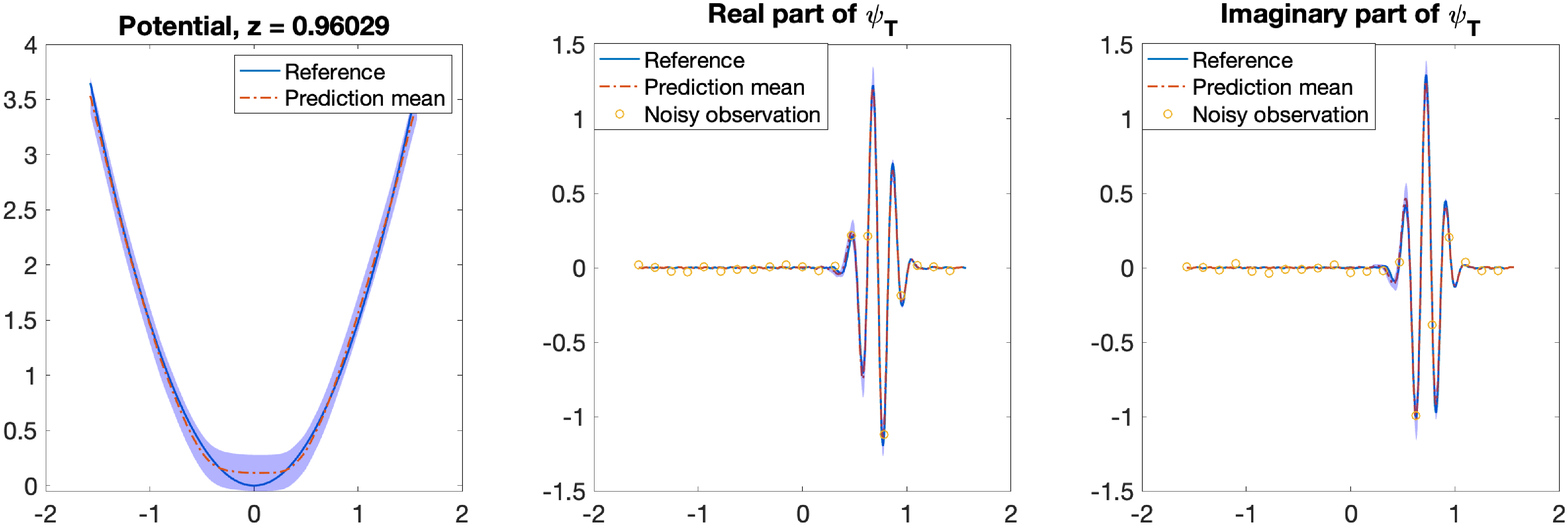}
	\caption{Test II, 20 sensors, true and predicted value of the potential function $\psi$ at final time $T=1.0$, for a training sample $z=0.0.7967$.}
	\label{fig:sto_psi1_50}
\end{figure}

\begin{figure}[!ht]
	\centering
	\includegraphics[width=0.9\textwidth]{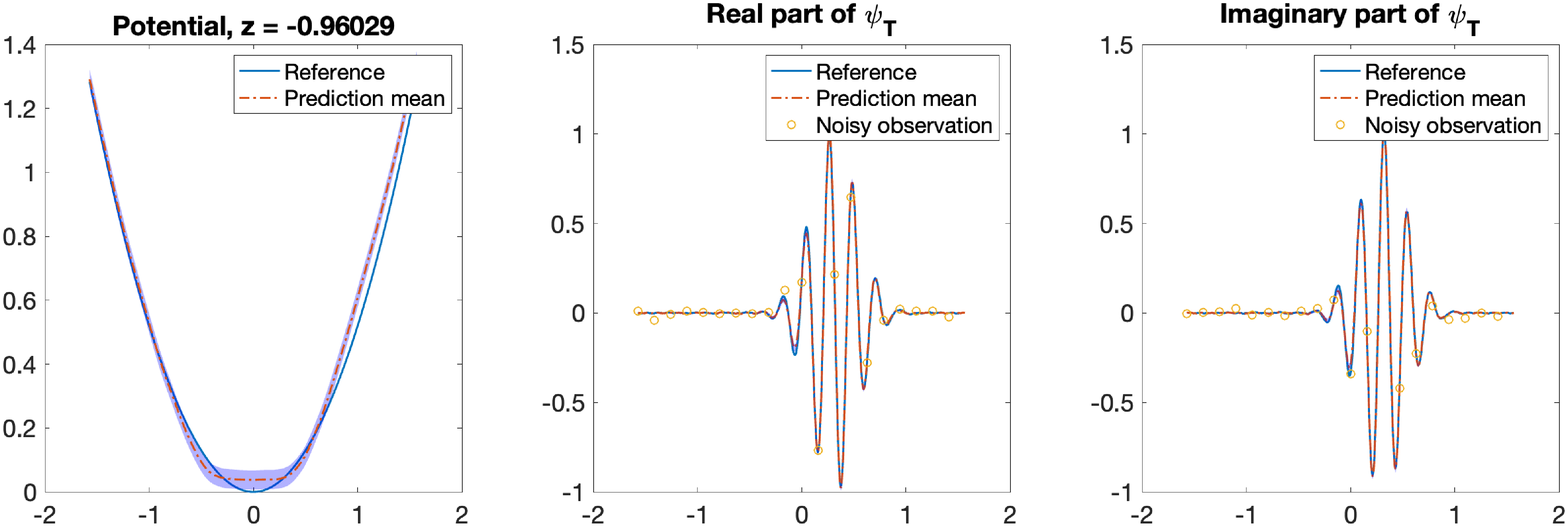}
	\caption{Test II, 20 sensors, true and predicted value of the potential function $\psi$ at final time $T=1.0$, for a training sample $z=-0.9603$.}
	\label{fig:sto_psi2_50}
\end{figure}

\begin{figure}[!ht]
	\centering
	\includegraphics[width=0.9\textwidth]{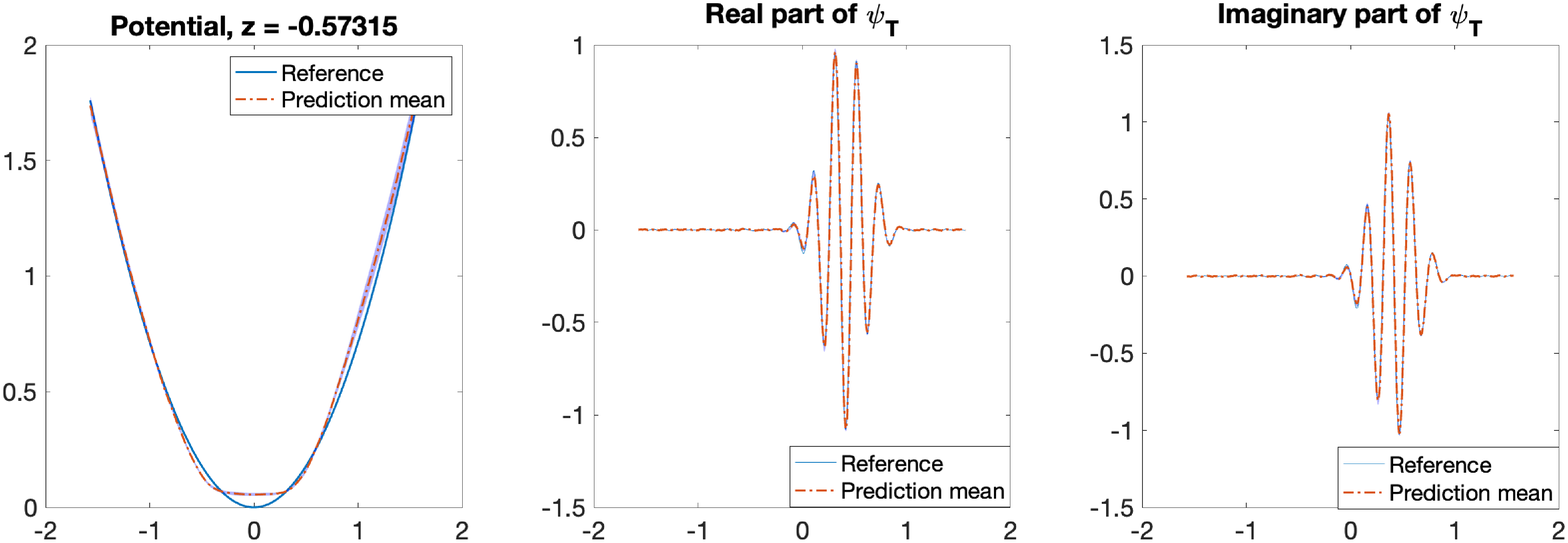}
	\caption{Test II, 50 sensors, testing case: $z=-0.57315$. True and predicted value of the potential function $\psi$ at final time $T=1.0$.}
	\label{fig:sto_psi_test_50}
\end{figure}

We then show the predictions of $V(x;z)$ when there are $20$ sensors and  $\psi_{\text{obs}}(x;z) \sim \mathcal{N}(\psi(x;z), 0.02)$, that is, the number of sensors are getting smaller and the noise in the observation data is also less. In this case, the predictions of the potential function for $z = [0.9603, 0.7967, 0.5255, 0.1834]$ are shown in Figure \ref{fig:sto_V_50}. The corresponding predictions of wave function $\psi$ with different values of $z = [0.9603, -0.9603]$ are shown in Figure \ref{fig:sto_psi1_50} and \ref{fig:sto_psi2_50}, respectively. In addition, a testing case for $z=-0.57315$ is presented in \ref{fig:sto_psi_test_50}. 
 
We observe that the results are still quite satisfactory under this test setting. This indicates our proposed network architecture and training algorithm can work well to learn the target stochastic potential, when the observation data is corrupted with a reasonable amount of noise.

\section{Conclusion and future work}

In this work, we adopt deep learning approach to learn the control variate in the inverse or control problem described by the Schr\"odinger equation in the semiclassical regime. With the choice of appropriate deep neural networks, we apply our framework to learn both the deterministic and stochastic control functions known as the potential. During the training process, the forward problem is solved by utilizing the efficient time-splitting spectral method, which guarantees the accuracy and enhances computational efficiency for the highly-oscillatory Schr\"odinger equation. We then develop a learning-based optimal control strategy by training neural networks to learn the control variate, considering observation data with or without noise. Our numerical results show that more reliable predictions can be obtained by adopting the SGLD algorithm.
{\color{black}We address the importance of our work by the following: (i) we investigate a {\it new} problem that is barely studied in the scientific computing fields; (ii) we introduce a {\it novel} hybrid NN-TSSP method as a deep learning approach to study the potential control problem described by the Schr\"odinger equation; (iii) the TSSP method as the forward solver in the sampling process is crucial, as the small parameter in the  Schr\"odinger equation brings numerical challenges. }

{\color{black}We mention some limitations of the current work, thus propose them as future works listed below.}
In the loss function during the training process, one can try to minimize the variance of the solution for more robust control. Besides, {\color{black}{we shall investigate higher-dimensional space problem for the Schr\"odinger equation, where other efficient schemes such as Gaussian wave packet based schemes can be adapted. 
}}Finally, more complicated potential function that depend on the temporal variable will be studied, in order to explore more general cases with practical applications for the quantum control problem.

\bibliographystyle{plain}
\bibliography{Ref}

\end{document}